\documentclass[11pt,openany,leqno]{article}
\usepackage[colorlinks=true]{hyperref}
\usepackage{amsmath,amsthm,amsfonts,amssymb,amscd,hyperref}
\usepackage[capitalize]{cleveref}
\usepackage{graphicx}
\usepackage{soul}
 
\numberwithin{equation}{section}
\input{xy} \xyoption{all}
\usepackage[utf8]{inputenc}
\usepackage{enumerate}

\usepackage{hyperref}
\usepackage{epsfig} 
\input{xy} \xyoption{all} 
\usepackage{mathrsfs}

\def \exp{\mathrm{exp}}

\def\qand{\quad \text{and}\quad}

\def\Ham{\mathrm{Ham}}
\def\Sym{\mathrm{Symp}}

\def\A{\mathbb A}

\def\E{\mathbb E}

\def\Q{\mathbb Q}

\def\C{\mathbb C}
\def\R{\mathbb R}
\def\N{\mathbb N}
\def\T{\mathbb T}
\def\D{\mathbb D}

\def\I{\mathbb I}
\def\I{\mathbb I}
\def\Z{\mathbb Z}
\def\cU{\mathcal U}

\def\det{\mathrm{det}}
\def\Re{\mathrm{Re}}

\def\leb{\mathrm{Leb}\, }
\def\Leb{\mathrm{Leb}\, }

\def\diam{\mathrm{diam}\:}

\def\se{{\mathsf e}}

\def\cal{\mathcal }
\def\sd{{\mathsf d}}

\newtheorem{proposition}{Proposition}[section]
\newtheorem{theorem}[proposition]{Theorem}
\newtheorem*{theorem*}{Theorem}
\newtheorem*{conjecture*}{Conjecture}
\newtheorem{corollary}[proposition]{Corollary}

\newtheorem*{problem*}{Problem}

\newtheorem{lemma}[proposition] {Lemma}

\newtheorem{theo}{Theorem}

\newtheorem{claim}[proposition]{Claim}
\newtheorem{fact}[proposition]{Fact}

\newtheorem{coro}[theo]{Corollary}

\newtheorem{question}[proposition]{Question}

\theoremstyle{remark}
\newtheorem{example}[proposition]{Example}
\newtheorem{remark}[proposition]{Remark}

\makeatletter

\DeclareTextFontCommand{\emph}{\em\bf}

\begin{document}
\title
{Analytic pseudo-rotations}
\author{Pierre Berger\thanks{IMJ-PRG, CNRS, Sorbonne Université, Université Paris Cité, partially supported by the ERC project 818737 Emergence of wild differentiable dynamical systems.
}}
\date{\today}
\maketitle
\begin{abstract} 
We construct analytic symplectomorphisms of the cylinder or the sphere with zero or exactly two periodic points that are not conjugate to a rotation. In the case of the cylinder, we show that these symplectomorphisms can be chosen to be ergodic or, to the contrary, with local emergence of maximal order. In particular, this disproves a conjecture of Birkhoff (1941) and solves a problem of Herman (1998). One aspect of the proof provides a new approximation theorem; it enables the implementation of the Anosov-Katok scheme in new analytic settings.
\end{abstract}

\tableofcontents

\section{Statement of the main theorems}

In low-dimensional analytic dynamics, a fundamental question is whether, far from periodic points, the dynamics is ``simple''.  In real or complex dimension one,  the dynamics restricted to an invariant domain without periodic points is either a rotation or in the basin of periodic points \cite{yoccoz1984n,sullivan1985quasiconformal}. In dimension 2, among volume dissipative complex or real analytic surface dynamics, the existence of  open subset  without periodic points but with very complicated  orbit has been recently revealed in \cite{BBiebler22}. For real analytic volume-preserving diffeomorphisms of the cylinder or the real sphere, 
a conjecture of Birkhoff \cite{Bi41}  states that such are also conjugate to rotations. We are going to disprove this conjecture. More precisely, we are going to show two examples of entire real symplectomorphisms of the annulus $\A= \R/\Z\times \R$, without periodic points and such that outside a region $\A_0'$ of $\A$ bounded by two disjoint analytic curves, the dynamics is analytically conjugate to a rotation, while inside $\A_0'$ the dynamics is not topologically conjugate to a rotation.
 In the first example, the restriction $f|\A_0'$ will be ergodic -- see \cref{main1} -- while in the second example it will be extremely far from being ergodic -- see \cref{main3}. More precisely, the ergodic decomposition of the latter map will be infinite-dimensional and even of maximal local order, i.e. with local emergence of maximal order~$2$. This confirms a conjecture on the typicality of high emergence in many categories \cite{Be_steklov_17}. 

From this we will deduce a corollary disproving the Birkhoff  conjecture. The proofs of the main theorems are shown by developing the Anosov-Katok method \cite{AK70}, together with a new approximation theorem for entire symplectomorphisms developing, a recent work with Turaev \cite{BT22}. The corollary regarding the sphere is obtained by using a blow-down technique introduced in \cite{berger2021coexistence}.

Finally, we will remark that these analytic constructions give examples of entire symplectomorphisms of $\C/\Z\times \C$ without periodic points and with a non-empty instability region $J$.

\subsection{Ergodic analytic and symplectic pseudo-rotation}
Let $S$ be an orientable analytic surface. We recall that a diffeomorphism of $S$ is \emph{symplectic} if it preserves the orientation and the volume. 
\begin{conjecture*}[Birkhoff {\cite[Pb 14-15]{Bi41}}]\label{Birkhoff}
An analytic symplectomorphism of the sphere with only two fixed points and no other periodic point is necessarily topologically conjugate to a rotation. 

Analogously, an analytic symplectomorphism of a compact cylinder without periodic points is topologically conjugate to a rotation.
 \end{conjecture*}

While following Birkhoff ``considerable evidence was adduced for this conjecture", Anosov-Katok \cite{AK70} gave examples of smooth symplectomorphisms of the sphere or the annulus with resp. 2 or 0 periodic points that are ergodic. Anosov and Katok proved their theorems by introducing the \emph{approximation by conjugacy method}, that we will recall in the sequel. Also, in his famous list of open problems, Herman wrote that a positive answer to the following question would disprove the Birkhoff  conjecture:
\begin{question}[Herman {\cite{He98}[Q.3.1]}]
Does there exist an analytic symplectomorphism of the cylinder or the sphere with a finite number of periodic points and a dense orbit?
\end{question}
\cref{maincoro} will disprove both Birkhoff's conjectures and answer positively to Herman's question in the cylinder case. It will also bring a new analytic case of application of the approximation by conjugacy method, as wondered by Fayad-Katok in \cite[\textsection 7]{FK04}. 
To state the main theorems, we~set:
\[ \T:= \R/ \Z ,\quad \I:= (-1,1) ,\quad \A:= \T\times \R, \qand \A_0= \T\times \I\; .\] 
We recall that a map of $\A$ is entire if it extends to an analytic map on $\C/\Z\times \C$. An \emph{entire symplectomorphism} is a symplectomorphism which is entire and whose inverse is entire.

An \emph{analytic cylinder} $\A_0'$ of $\A$ is a subset of the form $\{(\theta, y)\in \A: \gamma^- (\theta)< y< \gamma^+(\theta)\}$, for two analytic functions $\gamma^-<\gamma^+$. 

\begin{theo}\label{main1}
There is an entire symplectomorphism $F$ of $\A$ which leaves invariant an analytic cylinder $\A_0'\subset \A$, whose restriction to $\A_0'$ is ergodic and whose restriction to $\A\setminus cl(\A_0')$ is analytically conjugate to a rotation. Moreover, the complex extension of $F$ has no periodic points in $\C/\Z\times \C$. 
\end{theo}
 
The following consequence of \cref{main1} is the counterexample of both Birkhoff conjectures and a positive answer to Herman's question in the case of the cylinder:

\begin{coro}\label{maincoro}\begin{enumerate}
\item There is an analytic symplectomorphism of $cl(\A_0)$ which is ergodic and has no periodic points. 
\item There is an analytic symplectomorphism of the sphere, whose restriction to a sub-cylinder is ergodic and which displays only two periodic points.
\end{enumerate}\end{coro}

This corollary is proved in \cref{proof of corollaries} using the blow-down techniques introduced in \cite{berger2021coexistence}. 
 \cref{main1} will be proved by proving a new approximation \cref{approx analytic} which enables to implement the approximation by conjugacy method to analytic maps of the cylinder. Analytic symplectomorphisms of the 2-torus which are ergodic and without periodic points are known since Furstenberg \cite[Thm 2.1]{Fu61} (one of the explicit examples is $(\theta_1,\theta_2)\in \T^2\mapsto (\theta_1+\alpha,\theta_1+\theta_2)$ for any $\alpha$ irrational). Actually, the approximation by conjugacy method is known to provide examples of analytic symplectomorphisms of $\T^2$ which are isotopic to the identity, see \cite{BK19} for stronger results. Up to now, the only other known analytic realization of approximation by conjugacy method was Fayad-Katok's theorem \cite{FK14} showing the existence of analytic uniquely ergodic volume-preserving maps on odd spheres. As a matter of fact, we also answer a question of Fayad-Katok \cite[\textsection 7.1]{FK04} on whether analytic realization of the approximation by conjugacy method may be done on other manifolds than tori or odd spheres. In \cite{Be24}, we use the present work to enable various analytic and symplectic realizations of the approximation by conjugacy method on the sphere and the disk.

 A way to repair the Birkhoff conjecture might be to ask for the dynamics without periodic points\footnote{These are called pseudo-rotation in \cite{beguin2006pseudo}.} to have furthermore a rotation number satisfying a diophantine condition, as wondered by Herman in \cite{He98}[Q.3.2]. The approximation by conjugacy method has been useful to construct many other interesting examples of dynamical systems with special properties; see, for instance, the survey \cite{FH78,crovisier2006exotic}. Certainly the new approximation \cref{approx analytic} enables to adapt these examples to the case of analytic maps of the cylinder. In the next subsection,
we will state that this scheme enables to construct analytic symplectomorphisms with {maximal local emergence}:  an extreme  nonergodic behavior which has never been observed even in the  smooth case. 

\begin{center}
{\em The readers only interested in the Birkhoff conjecture can skip the next subsection and go directly to \cref{sketch}, and then skip \cref{sec:2.3,proof necklace}.}\end{center} 

\subsection{Analytic and symplectic pseudo-rotation with maximal local emergence} 
While an ergodic dynamics might sound complicated, the description of the statistical distribution of its orbits is by definition elementary. We recall that by Birkhoff's ergodic theorem, given a symplectic map $f$ of a compact surface $S$, for $\leb$ a.e. point $x$ the following limit is a well-defined probability measure called the \emph{empirical measure} of $x$.
\[\se(x):= \lim_{n\to \infty}\frac1{n} \sum_{k=1}^{n} \delta_{f^k(x) }\; .\]
The measure $\se(x) $ describes the statistical behavior of the orbit of $x$. The \emph{empirical function} $\se: x\in S\mapsto \se(x)$ is a function with value in the space $\cal M(S)$ of probability measures on $S$. Note that $\se$ is a measurable function. 

 A natural question is how complex is the diversity of statistical behaviors of the orbits of points in a Lebesgue full set. To study this, we shall look at the size of the pushforward $\se_*\leb$ of the Lebesgue probability measure $\leb$ of $S$ by~$\se$. The measure $\se_*\leb$ is called the \emph{ergodic decomposition}; it is a probability measure on the space $\cal M(S)$ of probability measures of~$S$. The ergodic decomposition describes the distribution of the statistical behaviors of the orbits. 
 
 To measure the size of the ergodic decomposition, for every compact metric space $X$, we endow the space $\cal M(X)$ of probability measures on $X$ with the \emph{Kantorovitch-Wasserstein metric $\sd$}:\label{defi Kantorovitch-Wasserstein}
 \[\forall \mu_1,\mu_2 \in \cal M(X), \quad \sd(\mu_1,\mu_2):= \inf\left\{ \int_{X^2} d(x_1,x_2) d\mu: \mu \in \cal M(X^2) \text{ s.t. }p_{i*} \mu=\mu_i\; \forall i\in \{1,2\}\right\},\]
 where $p_i: (x_1,x_2)\in X^2\to x_i\in X$ for $i\in \{1,2\}$. This distance induces the weak $\star$ topology on $\cal M(X)$ which is compact. Also it holds:
 \begin{proposition}\label{distance transport} For any compact metric spaces $X,Y$, any $\mu \in \cal M(X)$ and $f,g\in C^0(X,Y)$ it holds: \[\sd(f_*\mu,g_*\mu) \le \max_{x\in X} d(f(x),g(x)).\]
 \end{proposition}
\begin{proof} Let $\hat \mu$ be the measure on the diagonal of $X\times X$ which is pushed forward by the 1st and 2sd coordinate projections to $\mu$. Then observe that the pushforward $\nu$ of $\hat \mu$ by the product $(f,g)$ is a transport from $f_*\mu$ to $g_*\mu$; its cost $\int_{Y\times Y} d(y,y')d\nu $ is at most $\max_{x\in X} d(f(x),g(x))$. 
\end{proof}

We recall that that the \emph{$\epsilon$-covering number} of a compact metric space is the minimal number of $\epsilon$-balls needed to cover it. The following has been proved several times (see \cite[Thm 1.3]{BBochi21} for references):
 \begin{theorem} \label{boxorder} The metric order of $(\mathcal M(S),\sd) $ is $2$: 
 \[ \lim_{\epsilon\to 0} \frac{ \log \log \cal N(\epsilon)}{|\log \epsilon|}=2\quad \text{ with }\cal N(\epsilon)\text{ the $\epsilon$-covering number of }(\cal M(S),\sd)\; .\]
 \end{theorem}
In contrast, up to now, all the bounds on the dimensions of the ergodic decompositions of symplectic and analytic dynamics were finite: 
 \begin{example}[Case study]\label{casestudy}
 \begin{enumerate}
 \item If the measure $\leb$ is ergodic, then the ergodic decomposition is a Dirac measure at the Lebesgue measure: $\se_*\Leb= \delta_{\leb}$. 
 \item For a straight irrational rotation of the annulus $\A_0=\T\times \I$, we have $\se_*\Leb= \int_\I \delta_{\leb_{\T\times \{y\}} }d\leb(y)$ with $\Leb_{\T\times \{y\}}$ the one-dimensional Lebesgue measure on $\T\times \{y\}$. Hence the ergodic decomposition of a straight irrational rotation of the annulus is one-dimensional. The same occurs for an integrable twist map of $\A_0$. \end{enumerate}\end{example}
 
A natural problem (see \cite{Be_steklov_17,BeJMP22,BBochi21}) is to find dynamics for which $\se_*\Leb$ is infinite dimensional in many categories. The notion of emergence has been introduced to precisely state this problem. In the present conservative setting\footnote{The notion of emergence is also defined for dissipative system.}, the emergence describes the size of the ergodic decomposition. More precisely, the \emph{emergence} of $f$ at scale $\epsilon>0$ is the minimum number $\cal E(\epsilon)\ge 1$ of probability measures $(\mu_i)_{1\le i \le \cal E(\epsilon)}$ such that:
\[\int_{x\in S} \min_i d(\se(x), \mu_i)d\leb <\epsilon\; .\]
The \emph{order of the emergence} is 
\[\overline {\cal O}\cal E_f =\limsup_{\epsilon\to 0} \frac{ \log \log \cal E(\epsilon)}{|\log \epsilon|}\; . \]
In \cite[ineq. (2.2) and Prop. 3.14]{BBochi21}, it has been shown that $\overline {\cal O}\cal E_f \le \dim S=2$. Also, we constructed an example of $C^\infty$-flows on the disk with maximal emergence order: $ \overline {\cal O}\cal E_f =2$. 
We do not know how to perform this example in the analytic setting\footnote{Yet in the dissipative setting, a locally dense set of area-contracting polynomial automorphisms of $\R^2$ has been shown to have emergence of order 2 in \cite{BBiebler22}.}. Also, this example has an ergodic decomposition of local dimension 1. Actually, in view of \cref{boxorder}, one can hope for the existence of an ergodic decomposition of infinite local dimension and even a local emergence of positive order.
 
Let us denote $ \hat \se:= \se_*\leb$ the ergodic decomposition of $f$. The \emph{order of the local emergence} of $f$ is\footnote{The notion of local emergence appears first in \cite[Def 4.21 \& Pb 4.22]{BeJMP22}}:
\[ \overline{\cal O} \cal E_{loc}(f) = \limsup_{\epsilon \to 0} \int \frac{\log \left | \log \hat \se (B(\mu ,\epsilon) )\right| }{|\log \epsilon | } \, d\hat \se\; .\]
In \cite{He22}, Helfter showed that for $\hat \se$ a.e. $\mu\in \cal M(S)$ it holds:
\begin{equation}\label{Helfter_inequality}
\overline{\cal O}_{loc} \hat \se (\mu ):= \limsup_{\epsilon \to 0} \frac{\log \left | \log \hat \se (B(\mu ,\epsilon) )\right| }{|\log \epsilon | }\le \overline {\cal O}\cal E_f\; . 
\end{equation}
As $\overline {\cal O}\cal E_f \le 2$ by \cref{boxorder}, it comes that $\overline{\cal O}_{loc} \hat \se (\mu )\le 2$ a.e. Thus if $\overline{\cal O} \cal E_{loc}(f) =2$, then $\overline{\cal O}_{loc} \hat \se (\mu )= 2$ for $\hat \se$ almost every $\mu\in \cal M(S)$. 
In this work we give the first example of smooth symplectomorphism with infinite dimensional local emergence for smooth dynamics. Moreover, our example is entire and of maximal order of local emergence: 
\begin{theo}\label{main3}
There is an entire symplectomorphism $F$ of $\A$ which leaves invariant an analytic cylinder $\A_0'\subset \A$, whose restriction to $\A_0'$ 
has order of local emergence 2 and whose restriction to $\A\setminus cl(\A_0')$ is analytically conjugate to a rotation. Moreover, the complex extension of $F$ has no periodic points in $\C/\Z\times \C$. 
\end{theo}
A more restrictive version of local emergence could be done by replacing the $\liminf$ instead of $\limsup$. Then a natural open problem is:
\begin{question}Does there exist a smooth conservative map such that the following limit is positive:
\[\underline{\cal O} \cal E_{loc}(f) =\liminf_{\epsilon\to 0} \int \frac{\log \left | \log \hat \se (B(\mu ,\epsilon) )\right| }{|\log \epsilon | } \, d\hat \se\; ? \]
\end{question} 

\subsection{Sketch of proof and main approximation theorem} \label{sketch}
Let $\Sym^\infty(\A)$ denote the space of symplectomorphisms of $\A$. 
 Let $\Ham^\infty( \A)$ be the subgroup of Hamiltonian $C^\infty$-maps of $ \A$: there is a smooth family $(K_t)_{t\in [0,1]}$ of Hamiltonians $K_t\in C^\infty (\A)$ 
 which defines a family $(f_t)_{t\in [0,1]}$ such that $f_0=id$, $f_1=f$ and $\partial_t f_t$ is the symplectic gradient of $K_t$. Let $\Ham^\infty_0(\A)$ be the subgroup of $\Ham^\infty(\A)$ formed by maps 
 such that $(K_t)_{t\in [0,1]}$ can be chosen supported by $\A_0$. 
 We have $\Ham^\infty( \A)\subsetneq \Sym^\infty(\A)$. Indeed, the rotation $R_\alpha: (\theta,y)\in \A\mapsto (\theta+\alpha,y)$ of angle $\alpha\in \R$ belongs to $ \Sym^\infty(\A)\setminus \Ham^\infty( \A)$.

 The celebrated approximation by conjugacy method introduced by Anosov-Katok gives the existence of a map which satisfies all the properties of \cref{main1} but the analyticity. Our first challenge is to perform this construction among analytic maps of the cylinder. Let us recall:
\begin{theorem}[Anosov-Katok { \cite{AK70}}]\label{Anosov-Katok}
There exists a sequence $(H_n)_{n\in \N}$ of maps $H_n\in \Ham^\infty_0(\A)$ and a sequence of rational numbers $(\alpha_n)_{n\ge 0}$ converging to $\alpha\in \R\setminus \Q$ such that the sequence of maps $F_n := H_n\circ R_{\alpha_n}\circ H_n^{-1}$
converges to a map $F\in \Sym^\infty(\A)$ whose restriction $F|\A_0$ is ergodic. 
\end{theorem} 
We will give a complete proof of this theorem in \cref{proof of Anosov-Katok}. Let us here sketch its proof to understand the difficulty of adapting it to the analytic case (\cref{main1}). 
\begin{proof}[Sketch of proof] We construct $\alpha_n= p_n/q_n$ and $H_n$ by induction such that:
\begin{enumerate}[(1)]
\item for most points $y\in \I$, the pushforward of the Lebesgue measure on $\T\times \{y\}$ by $H_{n+1}$ is close to the Lebesgue measure on $\A_0$,
\item $H_{n+1}= H_{n} \circ h_{n+1}$ with $h_{n+1}$ which commutes with the rotation of angle $\alpha_n$. 
\end{enumerate}
The first property is obtained by constructing a map $h_{n+1}$ as depicted in \cref{figure_ergodic}, so that it sends most of the Lebesgue measure of $\T\times \{y\}$ nearby the Lebesgue measure of $\A_0$. 
\begin{figure}[h]
\centering
\includegraphics[width=13cm]{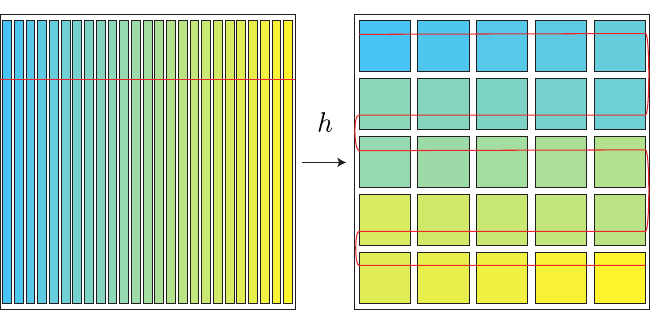}
\caption{A conjugacy which sends horizontal curves to curves close to 
be equidistributed with respect to the Lebesgue measure.}\label{figure_ergodic}
\end{figure}

As $H_n$ is area preserving and fixed before $h_{n+1}$, this property remains true for $H_{n+1}$ as claimed in (1). To obtain (2), it suffices to consider a lift of $h_{n+1}$ by the covering $(\theta, y)\mapsto (q_{n+1} \cdot \theta, y)$. 

Finally, we take $\alpha_{n+1}=p_{n+1}/q_{n+1}$ close to $\alpha_n$ so that $F_{n+1}$ is close to:
\[H_{n} \circ h_{n+1}\circ R_{\alpha_n} \circ h_{n+1} ^{-1} \circ H_{n} ^{-1}= H_{n} \circ R_{\alpha_n} \circ H_{n} ^{-1}=F_n\; ,\]
 as $h_{n+1}$ commutes with $R_{n+1}$. This implies the convergence of $(F_n)_n$ to a map $F$.

 Moreover, as $ q_{n+1}$ is large, all the points $(\theta,y)$ have their $R_{\alpha_{n+1}}$-orbits which are close to be equidistributed on $\T\times \{y\}$. Thus most of the points of $\A_0$ have their $F_{n+1}$-orbits close to be equidistributed on $\A_0$. This implies that $F_n$ is close to be ergodic and thus $F$ is ergodic. 
\end{proof} 

To disprove the Birkhoff conjecture using the approximation by conjugation method, we shall find a sequence of analytic maps $(h_n)_n$ satisfying $(1)$ and $(2)$. Actually the whole property $(1)$ is not needed, but we need at least that the sequence of compositions $H_n= h_1\circ h_2\cdots \circ h_n$ diverges on $\A_0$ without having a singularity in a uniform complex neighborhood of $\A_0\subset \C/\Z\times \C$ and likewise for its inverse. An idea to tackle this problem is to define $h_n$ as the time-one map of an entire Hamiltonian or to define it implicitly by an entire generating function. Then, in general, the map $h_n$ is not an entire automorphism: it or its inverse has a singularity in $\A_\C$. So we have to check that when $n$ is large,  the image by $H_n$ of this singularity does not approach $\A\subset \A_\C$. This is not easy. Step (ii) is even more problematic: after a $q_n$-lifting of $h_n$ to obtain $(ii)$, a singularity of $h_n$ becomes $1/q_n$-times closer to~$\A$. 

To avoid these problems, the idea is to work with \emph{entire real automorphisms}. Then there are no singularities to take care of. The complex definition domains of the maps and of their inverses cannot shrink since they are maximal, that is, equal to:
\[\A_\C:= \C/ \Z\times \C\; .\]
\emph{Horizontal twist maps} which are Hamiltonian maps of the form $(\theta, y)\mapsto (\theta+\tau(y), y)$ and \emph{vertical twist maps} which are Hamiltonian maps of the form $(\theta, y)\mapsto (\theta, y+v(\theta))$ provide examples of entire automorphisms. It suffices to take any $\tau $ and $v$ entire. Horizontal twist maps leave $\A_0$ invariant, while vertical twist maps do not. One can show that an entire symplectomorphism which leaves $\A_0$ invariant must be a horizontal twist map\footnote{ For every $y\in \A_0$, we apply Picard's Theorem to the second coordinate projection of the entire symplectomorphism restricted to $\C/Z\times \{y\}$ to conclude that it must be constant and so equal to $y$.}. To overcome the lack of such entire symplectomorphisms, we relax the condition on the invariance of $\A_0$ 
and we ask instead that the entire map should be close to the identity on a large compact subset of \[\E_\C:= \{(\theta,y) \in \C/\Z\times \C: \Re(y)\notin \I\}.\]
We will see below (see the main approximation \cref{approx analytic}) that there are a lot of such mappings. Then the idea is to use them in the approximation by conjugacy method (as the $h_n$). Indeed, the approximation by conjugacy method ensures the convergence of the $F_n= H_n \circ R_{\alpha_n} \circ H_n^{-1}$ on $\C/\Z\times \C$ to a map which is not conjugate to a rotation while we know that $(H_n)_n$ converges on a large compact subset of $\E_\C$ (a sequence of compositions of sufficiently close to the identity maps converges).

The main approximation \cref{approx analytic} is done by adapting to the analytic setting a new result with Turaev \cite{BT22}, which states that compositions of horizontal and vertical twist maps are dense in $\Ham^\infty(\A)$, see \cref{BT thm}. It is easy to approximate these maps by entire automorphisms, yet these might be far from preserving $\A_0$. That is why, we will show that the set of commutators of vertical and horizontal twist maps with the set of horizontal twist maps generates a dense subgroup of $\Ham^\infty(\A)$. Furthermore, we will notice that to approximate $\Ham^\infty_0(\A)$ only horizontal twist maps or such commutators supported by $\A_0$ are necessary. See \cref{approx reel1}. 

To pass to the analytic case, we will use Runge's theorem to approximate the latter maps by entire automorphisms which are arbitrarily close to the identity on large compact subsets of  $\A\setminus \A_0$. This leads to the main approximation \cref{approx analytic} stated below. This apparently simple step is actually one of the main technical difficulties of this work. \medskip 

We need a few notations to state \cref{approx analytic}. For every $\delta> 0$, let $\I_\delta:= [-1+\delta, 1-\delta]$, 
$ \A_\delta= \T\times \I_\delta$ and let $\Ham^\infty_\delta(\A)$ be the subgroup of $\Ham^\infty_0(\A)$ formed by flow maps 
of Hamiltonian $(H_t)_{t\in [0,1]}$ supported by $\A_\delta$. Observe that:\label{def hamdelta}
 \[\mathbb \A_0 =\bigcup_{\delta >0 } \A_\delta \qand \Ham^\infty_0(\A)=\bigcup_{\delta>0} \Ham^\infty_\delta(\A)\; .\] 

We denote $\Ham^\omega (\A) $ be the subgroup of $\Ham^\infty(\A)$ formed by Hamiltonian entire real automorphisms. For $\rho>1$, put:
\[K_\rho:= \T_\rho \times Q_\rho\quad \text{where }\T_\rho := \T+i [-\rho ,\rho]\qand Q_\rho:= [-\rho,-1 ]\sqcup [1 ,\rho] + i [-\rho ,\rho] 
\; .\]
We denote $\Ham_{\rho}^\omega (\A)$ the subset of $\Ham^\omega(\A)$ formed by maps whose restriction to $K_p$ is $\rho^{-1}$-close to the canonical inclusion $K_\rho\hookrightarrow \A_\C$:
 \[\Ham_{\rho}^\omega (\A):=\left\{h\in \Ham^\omega(\A): \; 
 \sup_{K_\rho } |h-id|<\rho^{-1}\right\}
\; .\]
\begin{theorem}[Main approximation result] \label{approx analytic}Let $0<\delta < 1 $, let $h\in \Ham_{\delta}^\infty(\A)$ and let $\cal U$ be a neighborhood of the restriction $h| \A_\delta $ in $C^\infty(\A_\delta, \A)$. 
Then for any $\rho> 1$, there exists $\tilde h\in \Ham_{\rho}^\omega (\A)$ such that the restriction $\tilde h | \A_\delta $ is in $\cal U$.
\end{theorem} 
\begin{remark} The proof of this theorem can be adapted to obtain an analogous statement for the surface $\T^2$ instead of $\A$.
\end{remark}
A consequence of \cref{approx analytic} regards $(1/q,0)$-\emph{periodic maps} for $q\ge 1$, which are maps $h$ satisfying $h(\theta+1/q,y)= h(\theta,y)+(1/q,0)$ for every $(\theta,y)\in \A$:
\begin{corollary} \label{coro:approx analytic}Let $0<\delta <1$, let $h\in \Ham_{\delta}^\infty(\A)$ be $(1/q,0)$-periodic and let $\cal U$ be a $C^\infty$-neighborhood of the restriction $h| \A_\delta $. 
Then for any $\rho>1$, there exists $\tilde h\in \Ham_{\rho }^\omega (\A)$ which is $(1/q,0)$-periodic and whose restriction $\tilde h | \A_\delta $ is in $\cal U$.
\end{corollary} 
\begin{proof} The map $h$ induces an entire automorphism $[h]$ on the quotient $\R/q^{-1} \Z\times \R$. Note that $\R/q^{-1} \Z\times \R$ is diffeomorphic to $\A$ via the map $\psi: (\theta, y)\in \R/q^{-1} \Z\times \R\to (q\cdot \theta, y)\in \A$ and that $\psi \circ [h]\circ \psi^{-1}$ satisfies the assumptions of \cref{approx analytic}. Hence there exists a map $\tilde g\in \Ham_{q\cdot \rho} (\A)$ such that the restrictions $\tilde g|\A_\delta$ and $\psi \circ [h]\circ \psi^{-1}$ are close. 
Then note that $\psi^{-1}\circ \tilde g\circ \psi $ is an entire automorphism on the quotient $\R/q^{-1} \Z\times \R$. It defines a map $\tilde h\in \Ham_\rho^\omega (\A)$ which is $1/q$-periodic and such that the restriction $\tilde h | \A_\delta $ is close to $h|\A_\delta$. \end{proof}
We will plug this corollary into the approximation by conjugation method to prove \cref{main1} in \textsection \ref{proof of main1}. 

\medskip 
To prove \cref{main3}, we will first implement the approximation by conjugation method to prove the following new result:
\begin{theorem}\label{main3 infty}
There is $F\in \Sym^\infty(\A_0)$ without periodic points and whose order of local emergence is 2. Moreover, $F$ is $C^0$-conjugate to a rotation. 
\end{theorem} 
We will give a complete proof of this theorem in \cref{proof of main3}. Let us sketch its proof.
\begin{proof}[Sketch of proof]
The proof follows the same lines as the one of Anosov-Katok, but differs at one point: we will not assume that for most $y\in \I$, the map $h_{n+1}$ pushes forward  the measure $\Leb_{\T\times \{y\}}$ to one near $\Leb_{\A_0}$. Instead, we will take the map $h_{n+1}$ so that for most of the point $y$, the measure of the set of $y'$ such that 
$h_{n+1*}\Leb_{\T\times \{y'\}}$ is $\epsilon_n$-close to $h_{n+1*}\Leb_{\T\times \{y\}}$ is smaller than $\exp(-\epsilon^{2-\eta_n}_n)$ for some $\epsilon_n,\eta_n\to 0$, see \cref{figure_emergence} and \cref{neck lace}. This lemma is proved using Moser's trick (see \cref{coro moser}) and a combinatorial \cref{coloriage} develloping \cite[Prop 4.2]{BBochi21}. \begin{figure}[h]
\centering
\includegraphics[width=13cm]{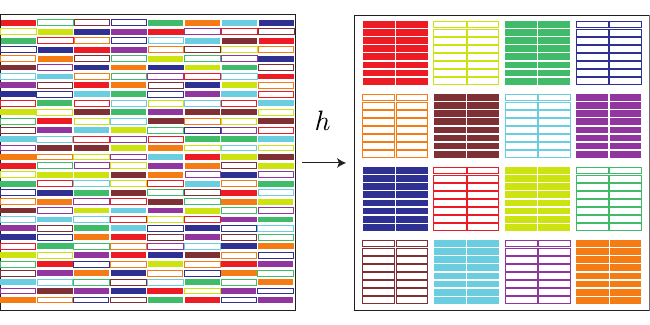}
\caption{A conjugacy which sends horizontal curves in different strips to distant ones in law.}\label{figure_emergence}
\end{figure}
\end{proof}
Again, we will plug \cref{coro:approx analytic} into this proof to show of \cref{main3} in \textsection \ref{proof of main3ana}. 
 \subsection{A complex analytic consequence of the main results}
Given a complex dynamics $f$ of several complex variables, there are two sets of special interest. 
The closure $J^*$ of the set of hyperbolic periodic points of $f$ and the set $J$ of points $z$ with bounded and Lyapunov unstable backward and forward orbits.
In the setting of polynomial automorphisms of $\C^2$, a natural  problem from the early 90s attributed to Hubbard
is whether $J=J^*$. In the conservative setting,
 this problem can be seen as a complex counterpart of Birkhoff's conjecture, since on each component of the $J$'s complement, the dynamics escapes to infinity or is analytically conjugate to a rotation, see \cite[Thm 2.1]{Be18}. 
In the dissipative setting, arguments in favor of $J=J^*$ are the recent results of Crovisier-Pujals \cite{CP18} and Dujardin \cite{Du20} which assert that for every real and resp. complex Hénon map satisfying $|\det D f|\cdot (\deg f)^2<1$, then 
 $J=J^*$ modulo a null measure set for any invariant probability measure. In the conservative setting, we observe that the complexified Furstenberg symplectomorphism $(\theta_1,\theta_2)\in \C^2/\Z^2
 \mapsto (\theta_1+\alpha,\theta_2+\theta_1)$ with $\alpha\in \R\setminus \Q$ satisfies $J= \R/\Z \times \C/\Z\neq \emptyset = J^*$. We bring the following new cases: 
 \begin{coro}\label{main2}The real entire symplectomorphisms $F$ of $\C/\Z\times \C$ of Theorems \ref{main1} and \ref{main3} have non-empty Julia sets $J$ which contain $\A_0'$ and has no  periodic points.
In particular, $J\neq J^*$ in each of these examples. 
\end{coro} This corollary is proved in \cref{proof of corollaries}. This leads to the following natural question:
\begin{question}Does there exist a conservative entire map of $\C^2$ which satisfies $J\neq J^*$ ?\end{question}
 
\medskip 
\emph{Acknowlegments.}
{\em 
I am grateful to Vincent Humili\`ere and Sobhan Seyfaddini for their advice on how to write the proof of Moser-Krygin's \cref{coro moser}. I am thankful to Romain Dujardin for several discussions around the $J=J^*$ problem. I am grateful to Bassam Fayad for his bibliographical and general advice. I am indebted to Raphael Krikorian for many comments on the redaction of the proof of the main theorem. I am thankful to the referees for all their valuable corrections and advice.}
 
 \section{Proof of the main theorems}

\subsection{Notations for probability measures} 
Let $(\cal M,\sd) $ denote the space of probability measures on $\A$ endowed with the Kantorovitch-Wasserstein distance (see \textsection \ref{defi Kantorovitch-Wasserstein}).
We denote by $\Leb$ the Lebesgue measure on $\A$ and by $\leb_\T$ the one-dimensional Lebesgue measure on $\T$.  
 For every cylinder $\A_0'\subset \A$, we set:
 \[\Leb_{\A'_0}:= \tfrac{1}{\leb \A_0'} \leb|\A'_0 ,\]
 which is a probability measure on $\A'_0$. For every $y\in \R$, we denote by $\leb_{\T\times \{y\}}$ the probability measure $\leb_{\T}\otimes \delta_{ y }$. For every $x\in \A$, $n\ge 1$ and $F\in \Sym^\infty(\A)$, 
 the \emph{$n^{th}$-empirical measure} is:
\[\se_{n}^F (x):= \frac1{n} \sum_{k=1}^{n} \delta_{F^k(x) }\; .\]
By Birkhoff's ergodic theorem, if a cylinder $\A_0'$ is left invariant by $F$ then for Leb. a.e. $x\in \A'_0$, the sequence $(\se^F_{n} (x))_n$ converges to a probability measure $\se^F(x)$. We recall that $F|\A'_0$ is ergodic iff $\se^F(x)=\Leb_{\A_0'}$ for a.e. $x\in \A_0'$.

\subsection{Proof of \cref{main1} on ergodic pseudo-rotations} 
To prove that \cref{approx analytic} implies the main theorems, we shall first recall the proof of Anosov-Katok's \cref{Anosov-Katok}.
\subsubsection{Existence of ergodic smooth pseudo-rotations}\label{proof of Anosov-Katok}
The proof of \cref{Anosov-Katok} is done by induction using a sequence $(\epsilon_n)_n$ of positive numbers satisfying: 
\[\sum_n \epsilon_n<1\; .\]
We recall that for $\alpha\in \R$, the rotation of angle $\alpha$ is denoted by $R_\alpha: (\theta,y)\in \A \mapsto (\theta+\alpha,y)$. Also for $\epsilon>0$, we recall that $\A_\epsilon:= \T\times [-1+\epsilon, 1-\epsilon]$. 
We prove below:
\begin{claim} \label{AK1}
There are sequences of rational numbers $\alpha_n= \tfrac{p_n}{q_n}$ and of Hamiltonian maps $H_n\in \Ham^\infty_0(\A)$ such that $F_n:= H_n\circ R_{\alpha_n} \circ H_n^{-1}$ satisfies for every $n\ge 1$:
\begin{enumerate}[$(P_1)$]
\item 
for every $x\in H_n(\A_{\epsilon_n})$, the empirical measure $\se^{F_n} (x)= \se^{F_n}_{q_n} (x)$ is $\epsilon_n$-close to $\leb_{\A_0}$. 
\item The $C^0$-distance between $\se^{F_n}_{k}$ and $\se^{F_{n-1}}_{k} $ is smaller than $\epsilon_n$ for every $k\le q_{n-1}$.
\item The $C^n$-distance between $F_n$ and $F_{n-1}$ is smaller than $\epsilon_n$:
\[ \max_{x\in \A, k\le n} \| D_x^k F_n -D_x^k F_{n-1} \|\le \epsilon_n\; .\]
\end{enumerate}
\end{claim} 
\begin{proof}[Proof that \cref{AK1} implies \cref{Anosov-Katok}]
 By $(P_3)$ and the convergence $\sum_n \epsilon_n$, the sequence $(F_n)_n$ converges to a  map $F\in \Sym^\infty(\A)$. 
 
Again, by convergence of $\sum\epsilon_n$ and the Borel-Canteli Lemma, the subset $B:=\bigcap_{n\ge 0} \bigcup_{k\ge n} H_k(\A_{\epsilon_k})$ of $\A_0$ has full Lebesgue measure in $\A_0$. Also by $(P_1)$ and $(P_2)$, for every $x\in B$, the measure $\leb_{\A_0}$ is an accumulation point of the sequence of empirical measures $(\se^{F}_n(x))_n$. Thus, by Birkhoff's ergodic theorem, a.e. point in $B $ satisfies that $\se^{F}_n(x)\to \Leb_{\A_0}$. As $B $ has full Lebesgue measure in $\A_0$ it comes that $F|\A_0$ is ergodic. \end{proof}
\begin{proof}[Proof of \cref{AK1}] Now let us show the claim by induction on $n$. Put $q_0=1$, $\alpha_0=0$ and $H_0=id$. 
Let $n\ge 1$, assume $(H_k)_{k<n} $ and $(\alpha_k=\frac{p_k}{q_k} )_{k<n}$ constructed and satisfying the induction hypothesis. 
We will show in \cref{proof Katok classic} the following:
\begin{lemma}\label{Katok classic}For every $\eta>0$, there exists a map $h_n\in \Ham^\infty_0(\A)$ such that:\begin{enumerate}[1.]
\item the map $h_n$ is $(1/q_{n-1},0)$-periodic. 
\item for every $(\theta,y)\in \A_{\epsilon_n}$, the map $ h_n$ pushes forward the one-dimensional Lebesgue probability measure $\Leb_{\T \times \{y\}}$ to a measure which is $\eta$-close to $\leb_{\A_0}$. 
\end{enumerate}
 \end{lemma}
 Then for $\eta$ sufficiently small, as $H_{n-1}$ leaves invariant $\A_0$, we obtain:
 \begin{enumerate}[\em 2'.] 
\item {\em for every $(\theta,y)\in \A_{\epsilon_n}$, the map $H_n:=H_{n-1}\circ h_n$ pushes forward the one-dimensional Lebesgue probability measure $\Leb_{\T \times \{y\}}$ to a measure which is $\epsilon_n/2$-close to $\leb_{\A_0}$.}
\end{enumerate}

Then by {\em 1.}, we have:
\[ H_n \circ R_{\alpha_{n-1}}\circ H_n^{-1} = H_{n-1} \circ h_n\circ R_{\alpha_{n-1}}\circ h_n^{-1}\circ H_{n-1}^{-1}=F_{n-1},\]
 and so for every $\alpha_n$ sufficiently close to $\alpha_{n-1}$, the map $F_n=H_n \circ R_{\alpha_{n}}\circ H_n^{-1} $ is sufficiently close to $F_{n-1}$ to satisfy $(P_2)$ and $(P_3)$. 
Also observe that when $\alpha_n=\frac{p_n}{q_n}$ is close to $\alpha_{n-1}$, then $q_n$ is large (with $q_n\wedge p_n=1$). For every $x= H_n(\theta,y)\in H_n (\A_{\epsilon_n})$, it holds that:
\[\se^{F_n}(x)= \frac1{q_n} \sum_{k=1}^{q_n} \delta_{H_n\circ R_{k\cdot p_n/q_n}(\theta,y)}=H_{n*} \left(\frac1{q_n} \sum_{k=1}^{q_n} \delta_{ R_{k\cdot p_n/q_n}(\theta,y)}\right)= H_{n*} \left(\frac1{q_n} \sum_{k=1}^{q_n} \delta_{(\theta+kp_n/q_n, y)}\right) \] 
is close to the pushforward by $H_n$ of $\Leb_{\T \times \{y\}}$ as $q_n$ is large. 
Thus by {\em 2'}, assumption $(P_1)$ holds true. 
\end{proof}
\begin{remark}\label{periodic point none}We can develop the construction to obtain immediately that $F$ has no periodic point. To this end, it suffices to add the following induction hypothesis\footnote{Actually 
using rotation number theory, one can show $F$ has no periodic point without this extra assumption, but this argument seems difficult to adapt to the complex setting.}: 
\begin{enumerate}[$(P_4)$]
\item for every $x\in \A_0$ and every $0<k< q_{n-1}$, it holds:
\[d( F^k_n(x), x)>d( F^k_{n-1}(x), x)\cdot (1-2^{-n})\; .\] 
\end{enumerate}
such an assumption is immediately verified by taking $\alpha_n$ sufficiently close to $\alpha_{n-1}$, because $d( F^k_{n-1}(x), x)$ is positive and $(1-2^{-n})<1$. 
 Also, it implies that for every $k\ge 0$, with $n$ such that $k< q_n$, it holds:
\[d( F^k(x), x)> d( F^k_{n}(x), x)\cdot \prod_{m\ge n} (1-2^{-m})>0\; .\] 
Consequently, $F$ has no periodic point.\end{remark}

\subsubsection{Existence of ergodic analytic pseudo-rotations} \label{proof of main1}
 We recall that $\E:= \A\setminus \A_0$ ,\quad $\E_\C:= \C/ \Z \times \{y\in \C: \Re(y)\notin \I\}$ and:
\begin{equation}\label{Krho recalled}K_\rho:= \T_\rho \times Q_\rho\quad \text{where }\T_\rho := \T+i [-\rho ,\rho]\qand Q_\rho:= [-\rho,-1 ]\sqcup [1 ,\rho] + i [-\rho ,\rho] 
\quad \forall \rho>1\; .\end{equation}

Similarly to the proof of \cref{Anosov-Katok}, we are going to show:
\begin{claim}\label{Claim CB1} Let $\rho>1$. 
 There exist  a sequence of entire automorphisms $H_n\in \Ham^\omega_\rho(\A)$ and a sequence of rational numbers $ \alpha_n$ such that:
\begin{enumerate}[$(i)$] 
\item The sequence $(\alpha_n)_n$ converges to a number $\alpha \in \R\setminus \Q$ and the sequence of restrictions $H_n|\E_\C$ converges to a map $H$ satisfying:
\begin{enumerate}
\item the restriction $H|int(\E_\C)$ is analytic,
\item the restriction $H|\C/\Z\times \{-1,1\}$ is analytic (as a one complex variable mapping),
\item it holds $\sup_{K_\rho} |H(x)-x|<\rho^{-1}$.
\end{enumerate} 
\item The sequence $(F_n)_n$ of maps $ F_n := H_n\circ R_{\alpha_n}\circ H_n^{-1}$ converges to an entire automorphism $F\in \Sym^\omega(\A)$ which is ergodic on $\A'_0:=\A\setminus H(\E)$. 
\item The map $F$ has no periodic points in $\C/\Z\times \C$. 
\end{enumerate}\end{claim}

\begin{proof}[Proof that \cref{Claim CB1} implies Theorems \ref{main1}]It remains only to show that $\A'_0:=\A\setminus H(\E)$ is an analytic cylinder. Note that the boundary of $cl(\A'_0)$ is equal to the one of $H(\E)$, and the latter is $H(\partial \E)= H(\T\times \{-1,1\})$. By (i), the restriction 
$H|\partial \E$ is analytic, and by the Cauchy inequality and $(i\cdot c)$, it is $C^1$-close to the canonical inclusion $\T\times \{-1,1\}\hookrightarrow \A$ when $\rho$ is large. Thus we can choose $\rho$ large enough so that $H(\partial \E)= H(\T\times \{-1,1\})$ is the union of two analytic graphs and so that 
$\A'_0$ is an analytic cylinder. 
\end{proof} 
Similarly to Anosov-Katok's proof, \cref{Claim CB1} is proved by induction using a sequence $(\epsilon_n)_n$ of positive numbers such that: 
\[\sum_n \epsilon_n<1\; .\]
For $\rho=n\ge 1$, \cref{Krho recalled} defines sets $\T_n$, $Q_n$ and $K_n= \T_n\times Q_n$. Put:
\[\D(n):= \{x\in \C: |x|<n\}.\]
 We are going to prove:
\begin{claim}\label{Claim CB2} For any $N\ge 1$, there are sequences of rational numbers $\alpha_n= \tfrac{p_n}{q_n}\to \alpha\in \R$ and of entire automorphisms $H_n\in \Ham^\omega(\A)$ such that with $\alpha_0=1$, $H_0=id$ and $F_n:= H_n\circ R_{\alpha_n} \circ H_n^{-1}$ for every $n\ge 1$, it satisfies:
\begin{enumerate}[$(P_0)$]
\item The restrictions $H_n|K_{n+N} $ and $H_{n-1}|K_{n+N} $ are $\epsilon_{n+N}$-$C^0$-close with $H_0=id$. 
\end{enumerate} 
\begin{enumerate}[$(P_1)$]
\item for every $x\in H_n(\A_{\epsilon_n} )$, the empirical measure $\se^{F_n} (x)= \se^{F_n}_{q_n} (x)$ is $\epsilon_n$-close to $ H_{n-1*}\leb_{\A_0}$. 
\item The distance between $\se^{F_n}_{k} (x)$ and $\se^{F_{n-1}}_{k} (x)$ is smaller than $\epsilon_n$ for every $k\le q_{n-1}$ and $x\in \T_n \times \D(n)$.
\item The $C^0$-distance between $F_n|\T_n \times \D(n)$ and $F_{n-1}|\T_n \times \D(n)$ is smaller than $\epsilon_n$. 
The $C^0$-distance between $F_n^{-1}|\T_n \times \D(n)$ and $F_{n-1}^{-1}|\T_n \times \D(n)$ is smaller than $\epsilon_n$. 
\item For every $x\in \T_n\times \D(n)$ and every $0<k< q_{n-1}$, it holds:
\[d( F^k_n(x), x)>d( F^k_{n-1}(x), x)\cdot (1-2^{-n})\; .\] 
\end{enumerate}\end{claim} 
\begin{proof}[Proof that \cref{Claim CB2} implies \cref{Claim CB1}.] 
Let $\rho>1$ be given by \cref{Claim CB1}. Let $N>\rho$ be large enough such that $\sum_{n\ge N} \epsilon_n<\rho^{-1}$.

Let us prove $(i)$. Observe that for every $x\in \E_\C$, there exists $N_0\ge N$ such that $x\in K_{N_0}$. As $K_{N_0} \subset K_{n}$ for every $n\ge N_0$, by $(P_0)$, the sequence $H_n(x)$ converges. Also, if $x$ is in $int(\E_\C)$ (resp. in $\C/\Z\times \{-1,1\}$), then there is a bidisk $D$ (resp. disk) containing $x$ and included in $K_{N_0+1}$ (resp. $K_{N_0+1}\cap \C/\Z\times \{-1,1\}$). So $(H_n(D ))_n$ is bounded by $(P_0)$. Thus $H$ is analytic on $int(\E_\C)$ and $\C/\Z\times \{-1,1\}$. Moreover, observe that:
\[
\sup_{x\in K_\rho} \|H(x)-x\|\le \sup_{x\in K_N} \|H(x)-x\| \le \sum_{n\ge 0} \sup_{x\in K_N} \|H_{n+1} (x)-H_n(x)\|\le \sum_{n\ge N} \epsilon_n \le \rho^{-1}\; .\]

Let us prove $(ii)$. First observe that for the same reason as in the proof that \cref{Claim CB1} implies \cref{main1}, for $N$ large enough, the set $H_n(\A_0)$ is an analytic cylinder for every $n$ and the sequence $(H_n(\A_0))_n$ converges to an analytic cylinder $\A_0':= \A\setminus H(\E)$. Thus ${H_n}_*\leb_{\A_0}$ converges to the normalized Lebesgue measure $\leb_{\A_0'}$ on the cylinder $\A'_0$ which equals $\bigcap_{n\ge 0} \bigcup_{k\ge n} H_k(\A_0)$: 
\[\A'_0=\bigcap_{n\ge 0} \bigcup_{k\ge n} H_k(\A_0) \qand {H_n}_*\leb_{\A_0}\rightharpoonup \leb_{\A_0'}\; .\]
By $(P_3)$, we observe\footnote{
The sequences $(F_n)_n$ and $(F_n^{-1})_n$ converge uniformly on each set $\T_k\times \D(k)$ to maps $F$ and $G$, which are therefore entire. Given any $z\in \C/\Z\times \C$, the point $F(z)$ belongs to a certain $ \T_k\times \D(k)$, and so does $F_n(z)$ for $n$ large enough. Then, by uniform convergence, we have $z= \lim_{n\to \infty} F^{-1}_n\circ F_n(z)= G\circ F(z)$. Hence $F=G^{-1}$ is indeed an entire automorphism. 
 } that $(F_n)_n$ converges to an entire real automophism $F$ from $\A$ onto $\A$. As its Jacobean is constantly equal to $1$, so it is actually an entire automorphism in $\Sym^\omega(\A)$. By $(P_0)$, the map $F$ is integrable on $H(int\, \E)$. By the Borel-Cantelli lemma, the subset $C:=\bigcap_{n\ge 0} \bigcup_{k\ge n} H_k(\A_{\epsilon_k})$ has full Lebesgue measure in $\A'_0=\bigcap_{n\ge 0} \bigcup_{k\ge n} H_k(\A_0)$. Also by $(P_1)$ and $(P_2)$, for every $x\in C$, the measure $\Leb_{\A_0'} $ is an accumulation point of the sequence of empirical measures $(\se^{F}_n(x))_n$. Thus a.e. point in $\A_0'$ satisfies that $\se^{F}_n(x)=\leb_{\A_0'}$ and so $F|\A_0'$ is ergodic. 

Finally, using the same argument as in \cref{periodic point none}, we obtain that $(P_4)$ implies $(iii)$. 
\end{proof} 
 
 \begin{proof}[Proof of \cref{Claim CB2}] Let us show the claim by induction on $n\ge 0$. Put $q_0=1$, $\alpha_0=0$ and $H_0=id$. 
Let $n\ge 1$ and assume $(H_k)_{k<n} $ and $(\alpha_k=\frac{p_k}{q_k} )_{k<n}$ constructed and satisfying the induction hypothesis. 
By \cref{Katok classic}, there exists $\mathring h_n\in \Ham^\infty_0(\A)$ such that:\begin{enumerate}[(\em 1)]
\item the map $\mathring h_n$ is $(1/q_{n-1},0)$-periodic. 
\item for every $(\theta,y)\in \A_{\epsilon_n}$, the map $\mathring H_n:=H_{n-1}\circ \mathring h_n$ pushes forward the one-dimensional Lebesgue probability measure $\Leb_{\T \times \{y\}}$ to a measure which is $\epsilon_n/2$-close to $H_{n-1*}\leb_{ \A_0 }$. 
\end{enumerate}
Let $\delta>0$ be such that the support of $\mathring h_n$ is included in $\A_\delta$. 
By \cref{coro:approx analytic}, for any $C^\infty$-neighborhood $\cal U$ of the restriction $\mathring h_n | \A_\delta $, for any $\rho>1$, there exists $h_n\in \Ham_\rho^\omega (\A)$ which is $(1/q_n,0)$-periodic and such that the restriction $ h_n | \A_\delta $ is in $\cal U$. Now put:
\[H_n:= H_{n-1}\circ h_n\; .\]
Observe that if we fix $\rho>1$ large enough depending on $H_{n-1}$ and $\epsilon_{n+N}$, then the restrictions $H_n|K_{n+N} $ and $H_{n-1}|K_{n+N} $ are $\epsilon_{n+N}$-close, as stated by 
$(P_0)$. Furthermore, the measure $H_{n*} \Leb_{\A_0}=\Leb_{H_n(\A_0)}$ is $\epsilon_n/4$-close to $H_{n-1*}\Leb_{\A_0}$. We will fix $\cU$ in the sequel. Again, by $(1/q_n,0)$-periodicity of $h_n$, we have:
\[ H_n \circ R_{\alpha_{n-1}}\circ H_n^{-1} = H_{n-1} \circ h_n\circ R_{\alpha_{n-1}}\circ h_n^{-1} \circ H_{n-1}^{-1}=F_{n-1},\]
 and so for every $\alpha_n$ sufficiently close to $\alpha_{n-1}$, $F_n=H_n \circ R_{\alpha_{n}}\circ H_n^{-1} $ is sufficiently close to $F_{n-1}$ to satisfy $(P_2)$, $(P_3)$ and $(P_4)$. 
Also observe that when $\alpha_n=\frac{p_n}{q_n}$ is close to $\alpha_{n-1}$, then $q_n$ is large (with $q_n\wedge p_n=1$). Thus, for every $x= H_n(\theta,y)\in H_n (\A_{\epsilon_n})$, it holds that
\[\se^{F_n}(x)= \frac1{q_n} \sum_{k=1}^{q_n} \delta_{H_n\circ R_{k\cdot p_n/q_n}(\theta,y)}=H_{n*} \left(\frac1{q_n} \sum_{k=1}^{q_n} \delta_{ R_{k\cdot p_n/q_n}(\theta,y)}\right)= H_{n*} \left(\frac1{q_n} \sum_{k=1}^{q_n} \delta_{(\theta+kp_n/q_n, y)}\right) \] 
is $\epsilon_n/8$-close to the pushforward by $H_n$ of $\Leb_{\T \times \{y\}}$. 
Then assuming $\cU$ small enough, we obtain that $\se^{F_n}(x)$ is $\epsilon_n/4$-close to the pushforward by $\mathring H_n$ of the one-dimensional Lebesgue probability measure $\Leb|\T \times \{y\}$. Then by {\em (2)}, $\se^{F_n}(x)$ is $\tfrac34 \epsilon_n$-close to $H_{n-1*}\leb_{ \A_0 }$. As $H_{n*}\leb_{ \A_0 }$ and $H_{n-1*}\leb_{\A_0}$ are $\epsilon_n/4$-close, Property $(P_1)$ holds true.
\end{proof}
\subsection{Proof of \cref{main3} on pseudo-rotations with maximal local emergence}\label{sec:2.3}

The proof of \cref{main3} follows the same lines as the one of \cref{main1}. We shall first give the proof of \cref{main3 infty}, which is the smooth version of \cref{main3}. 
\subsubsection{Existence of smooth pseudo-rotations with maximal local emergence} \label{proof of main3}
The proof of \cref{main3 infty} is also done by induction using a sequence $(\epsilon_n)_n$ of positive numbers satisying: 
\[\sum_n \epsilon_n<1\; .\]
We prove below:
\begin{claim} \label{Claim E0}
There are sequences of rational numbers $\alpha_n= \tfrac{p_n}{q_n}$ converging to $\alpha\in \R $, a decreasing sequence of positive numbers $(\eta_n)_n$ converging to $0$ and a sequence of Hamiltonian maps $H_n\in \Ham^\infty_0(\A)$ such that $H_0=id$, $\alpha_0=1$ and for $n\ge 1$ it holds:
\begin{enumerate}[$(P_1)$]
\item There is a family $(I_{n,i})_{1\le i\le M_n} $ of $M_n\ge \exp( \eta_n^{-2+\epsilon_n}) $ disjoint segments of $\I$ s.t.:
\begin{enumerate}[(a)]
\item for any $i$,
 it holds $(1-\frac{\epsilon_n}{ M_n}) \frac{2 }{M_n}< \leb I_{n,i}< \frac{2 }{M_n} $; hence $\leb(\I\setminus \bigcup_{i} I_{n,i})<\frac{2}{M_n}\epsilon_n $,
 \item for any $y\in I_{n,i}$ and $y'\in I_{n,j}$, $j\neq i$, the distance between 
 $H_{n*} \leb_{\T\times \{y \} }$ and $H_{n*} \leb_{\T\times \{y'\}}$ is $> 2\eta_n $. 
 \end{enumerate} 
\item It holds:
\[
\sup_{x\in \A} \|H_n(x)-H_{n-1}(x)\|< \epsilon_{n-1}\eta_{n-1} \qand \sup_{x\in \A} \|H_n^{-1}(x)-H_{n-1}^{-1}(x)\|< \epsilon_{n-1}\eta_{n-1}\; . \]
\item The $C^n$-distance between $F_n:= H_n\circ R_{\alpha_n} \circ H_n^{-1}$ and $F_{n-1}$ is smaller than $\epsilon_n$. 
\end{enumerate}
\end{claim} 
 
\begin{proof}[Proof that \cref{Claim E0} implies \cref{main3 infty}] Let $(F_n)_n$ be given by the claim. By $(P_3)$ and the convergence $\sum_n \epsilon_n$, the sequence $(F_n)_n$ converges to a $C^\infty$-map $F\in \Sym^\infty(\A)$. Let us show that $F$ satisfies the properties of \cref{main3 infty}.

By $(P_2)$, the sequence $(H_m)_m$ converges to a homeomorphism $H$. As each $H_n$ coincides with the identity on the complement of $\A_0$, so does $H$ and so $F$ coincides with $R_\alpha$ on $\A\setminus \A_0$. 
Note also that $H$ preserves the volume. Also, by continuity of the composition, it holds:
\[F=H\circ R_\alpha\circ H^{-1},\] 
 and so $F$ is indeed $C^0$-conjugate to an irrational rotation. 
 
 Now let us show that $F$ has maximal local emergence. For every $n,i$, let $B_{n,i} := \T\times I_{n,i}$ and $B_n:=\bigcup_i B_{n,i}$. The set $C:=\bigcap_{n\ge 0} \bigcup_{k\ge n} H_k(B_k)$ is included in $\A_0$. By $(P_1)-(a)$, the measure of $ \bigcup_{k\ge n} H_k(B_k)$ is equal to the one of $\A_0$ and so 
 the subset $C $ has full Lebesgue measure in $ \A_0 $. 
 
 Now let $x \in C$. Then there exist $n$ arbitrarily large and $i$ such that $x$ belongs to $H_n(B_{n,i} )$. Put $x= H_n(\theta,y)$. Then by $(P_1)$, we have:
 \[\leb \left\{ y'\in \I : d(H_{n*}( \leb_{\T\times \{y\}}), H_{n*}( \leb_{\T\times \{y'\}}) )<2\eta_n \right\} 
< \leb I_{n,i} + \leb (\I\setminus \cup_j I_{n,j} )\le \frac{2(1+\epsilon_n)}{M_n} \; .\]
 By $(P_2)$ and \cref{distance transport}, for every $m>n$, it holds:
\[\leb \left\{ y'\in \I : d(H_{m*} ( \leb_{\T\times \{y\}}), H_{m*}( \leb_{\T\times \{y'\}}) )<(2-\sum_{k= {n}}^{m-1} \epsilon_k) \eta_n \right\} 
< \frac{ 2(1+\epsilon_n) }{M_n}\; .\]
 Taking the limit as $m\to \infty$ and inferring that $\sum_{k } \epsilon_k<1$, it comes:
\[\leb \left\{ y'\in \I : d(H_{*} ( \leb_{\T\times \{y\}} ), H_{*}( \leb_{\T\times \{y'\}} ) )< \eta_n \right\} 
< \frac{ 2 (1+\epsilon_n) }{M_n}\; .\]
 
As $F= H\circ R_\alpha\circ H^{-1}$ with $\alpha\in \R\setminus \Q$, the empirical measure $\se^F(x')$ of $x'= H(\theta',y')\in \A_0$ is equal to $H_* ( \leb_{\T\times \{y'\}})$. Thus:
\[\leb \left\{ (\theta',y')\in \A_0 : d(\se^F\circ H(\theta,y), \se^F\circ H(\theta',y') )< \eta_n \right\} 
< \frac{2 (1+\epsilon_n) }{M_n}\; .\]
Now we infer that $x=H(\theta,y)$ and that $H$ is area preserving to obtain:
\[\leb \left\{ x'\in \A_0 : d(\se^F (x), \se^F (x') )< \eta_n \right\} 
< \frac{2 (1+\epsilon_n) }{M_n}\; .\]

Thus
\[\frac{\log |\log \leb \{ x'\in \A_0 : d(\se ^{F} (x) , \se ^{F} (x') )\le \eta_n \}| }{-\log \eta_n} \ge 
\frac{\log | \log (2(1+\epsilon_n) ) - \log M_n| }{-\log \eta_n}\; .\]
As $M_n>\exp( \eta_n^{-2+\epsilon_n}) $, we obtain:
\[\frac{\log |\log \leb \{ x'\in \A_0 : d(\se ^{F} (x) , \se ^{F} (x') )\le \eta_n \}| }{-\log \eta_n} \ge \frac{\log ( \eta_n^{-2+\epsilon_n} - \log (2(1+\epsilon_n) )) }{-\log \eta_n} \sim 2
\; .\]
 So the local order is a.e. $\ge 2$. By Hefter's inequality \eqref{Helfter_inequality} \cpageref{Helfter_inequality} and \cref{boxorder}, the local order is at most $2$ and so equal to $2$.
 \end{proof}

 \begin{proof}[Proof of \cref{Claim E0}] 
The proof is done by induction on $n\ge 0$. The step $n=0$ is obvious. Assume $n\ge 1$, with $H_{n-1}$ and $\alpha_{n-1} =\frac{p_{n-1}}{q_{n-1} }$ constructed. 
 
 We are going to construct $H_n$ of the form $H_n=H_{n-1}\circ h_n$ with $h_n\in \Ham_0^\infty(\A)$ which is $(1/q_{n-1},0)$-periodic. Then, for $\alpha_n$ sufficiently close $\alpha_{n-1}$ we observe that $(P_3)$ is satisfied by Anosov-Katok's trick: for the $C^0$-compact open topology of $\C/\Z\times \C$ when $\alpha_n$ is close to $\alpha_{n-1}$ at $H_n$ fixed, it holds:
 \[ F_n =H_n\circ R_{\alpha_n} \circ H_n^{-1}\approx H_n\circ R_{\alpha_{n-1}} \circ H_n^{-1}= H_{n-1}\circ R_{\alpha_{n-1}} \circ H_{n-1}^{-1}=F_{n-1}\; .\]
 
 Consequently, to prove the claim, it suffices to find a $(1/q_{n-1},0)$-periodic $h_n\in \Ham_0^\infty(\A)$ such that $H_n=H_{n-1}\circ h_n$ satisfies $(P_1)$ and $(P_2)$. 
 
 Let $Q$ be a large multiple of $q_{n-1}$ and let $( T_{k,\ell})_{0\le k,\ell<Q}$ be the tiling of $\A_0$ defined by:
 \[ T_{k,\ell} := \left\{\theta\in \T : \tfrac k Q< \theta < \tfrac{k+1}Q \right\}\times \left(-1+\tfrac {2\ell } Q,-1+\tfrac{2\ell+2}Q\right)\; .\] 
When $Q$ is large and if $h_n\in \Ham_0^\infty(\A)$ leaves invariant the tiling, i.e. $h_n(T_{k,\ell})=T_{k,\ell}$ for every $(k,\ell)$, then $h_n$ and $h_n^{-1}$ are $Q^{-1}$-$C^0$ close to the identity. We fix $Q$ sufficiently large such that $H_n=H_{n-1}\circ h_n$ and its inverse are $\epsilon_{n-1}\eta_{n-1}$-$C^0$-close to respectively $H_{n-1}$ and its inverse, as claimed in $(P_2)$.
Consequently, to prove the claim, it suffices to find a $(1/q_{n-1},0)$-periodic $h_n\in \Ham_0^\infty(\A)$ which leaves invariant 
$( T_{k,\ell})_{0\le k,\ell<Q}$ and such that $(P_1)$ is satisfied. In order to do so, we use the following Lemma is shown in \cref{proof necklace}:
 \begin{lemma} \label{neck lace}
 For every $1>\epsilon>0$, there exist $\eta>0$ arbitrarily small, an integer $M\ge \exp( \eta^{-2+\epsilon})$, a map $h\in \Ham^\infty(\T\times \R/2\Z)$ and a family $( J_{\epsilon,i})_{1\le i\le M } $ of disjoint compact segments of $\R/2\Z$ such that:
\begin{enumerate}[(a)]
\item it holds $\leb(\I\setminus \bigcup_i J_{\epsilon,i} )=\frac{2}{M}\epsilon $ and $ \leb J_{\epsilon,i}=(1-\frac{\epsilon}{ M} ) \frac{2 }{M} $ for every $i$,
 \item for any $j\neq i$, $y\in J_{\epsilon,i} $ and $y'\in J_{\epsilon,j}$, the distance between the measures $h_* \leb_{\T\times \{y\}} $ and $h_* \leb_{\T\times \{y'\}}$ is greater than $2\eta$,
 \item the map $h$ coincides with the identity on neighborhoods of $\{0\}\times \R/2\Z$ and $\T\times \{0\}$. 
 \end{enumerate} 
 \end{lemma}

 Let us apply \cref{neck lace} with $\epsilon=\epsilon_n/2$ fixed and some $ \eta\in (0,\eta_{n-1})$ small enough so that 
\begin{equation}\label{minorationM}M\ge \exp( \eta^{-2+\epsilon_n/2})\ge \exp( \eta_n^{-2+\epsilon_n})\, \quad \text{ with } \eta_n:= \eta/(Q\cdot \|DH_{n-1}\|_{C^0}) \; . \end{equation}This is possible since the powers $-2+\epsilon_n/2$ and $-2+\epsilon_n$ are different in the latter inequality.

 Let $( J_{\epsilon,i})_{1\le i\le M} $ be the family of disjoint segments of $\I$ and let $h$ be the map provided by the lemma. For every $0\le k,\ell<Q$, let:
 \[ \phi_{k,\ell}: (\theta, y)\in T_{k,\ell}\mapsto (Q\cdot \theta, Q \cdot y+Q-2\ell)\in \T\times \R/2\Z\; . \]
Note that $\phi_{k,\ell}$ is a bijection from $T_{k,\ell}$ onto $(\T\setminus\{0\}) \times (\R/2\Z \setminus \{0\})$. We define:
 \begin{equation}\label{defhn}h_n:= (\theta,y)\mapsto \begin{cases} 
 (\theta,y) & \text{if } (\theta,y)\notin \bigcup_{0\le k,\ell<Q } T_{k,\ell}\, ,\\
 \phi_{k,\ell}^{-1} \circ h\circ \phi_{k,\ell}(\theta,y) & \text{if } (\theta,y)\in T_{k,\ell} \; .
 \end{cases} \end{equation}
 Remark that $h_n$ leaves invariant the tiling $( T_{i,j})_{0\le i,j<Q}$. Now let us verify $(P_1)$. First let $M_n:= Q\cdot M$. For $1\le i\le M_n$, put $i=:mM+\ell$ with $1\le \ell \le M$ and $0\le m<Q$. By identifying $J_{\epsilon, \ell}$ to a subset of $(0,2)$ we define:
 \[I_{n,i} := J_{\epsilon, \ell}/Q -1+2 m/Q \subset (-1,1)=\I\; .\]
First observe that by \cref{minorationM}, it holds:
\[M_n := Q\cdot M\ge Q\cdot \exp( \eta^{-2+\epsilon_n/2})\ge \exp( \eta_n^{-2+\epsilon_n}),\]
as claimed in the begining of $(P_1)$. Moreover by \cref{neck lace} (a), Property $(P_1)-(a)$ is verified as $\epsilon_n\ge \epsilon$. Also, for any $y\in I_{\epsilon,i}\neq I_{\epsilon,j}\ni y'$, the distance between the measures $h_{n*} \leb_{\T\times \{y\}} $ and $h_{n*} \leb_{\T\times \{y'\}}$ is greater than $2\eta/Q$ and so the distance between the measures $H_{n*} \leb_{\T\times \{y\}} $ and $H_{n*} \leb_{\T\times \{y'\}}$ is greater than $2\eta/(Q\cdot \|DH_{n-1}\|_{C^0})=2\eta_n $ as claimed in $(P_1)-(b)$. 
 \end{proof} 
\subsubsection{Existence of analytic pseudo-rotations with maximal local emergence} \label{proof of main3ana}
For the same reasons, \cref{main3} is an immediate consequence of the following counterpart of \cref{Claim CB1} where only $(ii)$ has been changed. 
\begin{claim}\label{Claim E1} Let $\rho>1$. 
 There exist  a sequence of entire automorphisms $H_n\in \Ham^\omega(\A)$ and a sequence of rational numbers $ \alpha_n$ such that:
\begin{enumerate}[$(i)$] 
\item The sequence $(\alpha_n)_n$ converges to a number $\alpha \in \R\setminus \Q$ and the sequence of restrictions $H_n|\E_\C$ converges to a map $H$ which is analytic on $int(\E_\C)$ and on $\C/\Z\times \{-1,1\}$ with $\sup_{K_\rho} |H(x)-x|<\rho^{-1}$.
\item The sequence $(F_n)_n$ of maps $ F_n := H_n\circ R_{\alpha_n}\circ H_n^{-1}$ converges to an entire automorphism $F\in \Sym^\omega(\A)$ whose restriction to $\A':=\A\setminus H(\E)$ has local emergence $2$. 
\item The map $F$ has no periodic points in $\C/\Z\times \C$. 
\end{enumerate}\end{claim}
Similarly, \cref{Claim E1} is proved below using some series $\sum_n \epsilon_n<1$ of positive numbers and the following counterpart of \cref{Claim CB2}:
\begin{claim}\label{Claim E2} For any $N\ge 1$, there are sequences of rational numbers $\alpha_n= \tfrac{p_n}{q_n}\to \alpha\in \R\setminus \Q$, of positive numbers $(\eta_n)_n$ converging to $0$ and of entire automorphisms $H_n\in \Ham^\omega(\A)$ such that $\alpha_0=1$, $H_0=id$ and $F_n:= H_n\circ R_{\alpha_n} \circ H_n^{-1}$, it holds for every $n\ge 1$:
\begin{enumerate}[$(P_0)$]
\item The restrictions $H_n|K_{n+N} $ and $H_{n-1}|K_{n+N} $ are $\epsilon_{n+N}$-$C^0$-close, with $H_0=id$. 
\end{enumerate} 
\begin{enumerate}[$(P_1)$]
\item There is a family $(I_{n,i})_{1\le i\le M_n} $ of $M_n\ge \exp( \eta_n^{-2+\epsilon_n}) $ disjoint subsegments of $ [-1+\delta_n ,1-\delta_n]$, with $\delta_n= {\tfrac12 \epsilon_n/M_n}$ s.t.:
\begin{enumerate}[(a)]
\item for any $i$,
 it holds $(1-\epsilon_n) \frac{2 }{M_n}< \leb I_{n,i}< \frac{2 }{M_n} $; hence $\leb(\I\setminus \bigcup_{i} I_{n,i})<\frac{2}{M_n}\epsilon_n $,
 \item for any $y\in I_{n,i}$ and $y'\in I_{n,j}$, $j\neq i$, the distance between 
 $H_{n*} \leb_{\T\times \{y \} }$ and $H_{n*} \leb_{\T\times \{y'\}}$ is $> \eta_n $. \end{enumerate} 
 \item It holds $\sup_{x\in \A_{\epsilon_n/M_n}} \|H_n(x)-H_{n-1}(x)\|< \epsilon_{n-1}\eta_{n-1}$.
\item The $C^0$-distance between $F_n|\T_n \times \D(n)$ and $F_{n-1}|\T_n \times \D(n)$ is smaller than $\epsilon_n$. 
The $C^0$-distance between $F_n^{-1}|\T_n \times \D(n)$ and $F_{n-1}^{-1}|\T_n \times \D(n)$ is smaller than $\epsilon_n$. 
\item For every $x\in \T_n\times \D(n)$ and every $0<k< q_{n-1}$, it holds:
\[d( F^k_n(x), x)>(1-2^{-n})\cdot d( F^k_{n-1}(x), x)\; .\] 
\end{enumerate}\end{claim} 
\begin{proof}[Proof that \cref{Claim E2} implies \cref{Claim E1}]
 Properties $(i)$ and $(iii)$ are proved as in the proof  that \cref{Claim CB2} implies \cref{Claim CB1}, and likewise for the fact that $(F_n)_n$ converges to an entire automorphism $F\in \Sym^\omega(\A)$. 
Using verbatim the same argument as in the proof of \cref{Claim E0} we obtain that the local emergence of $F|\A_0'$ is $2$ as stated in $(ii)$. 
\end{proof} 
\begin{remark} 
Note that $(P_2)$ implies that $(H_n|\A_0)_n$ converges for the $C^0$-compact-open topology. With slightly more work, one can also obtain that this limit is a homeomorphism onto its image. However, the construction does not imply that this limit extends $H: \E\hookrightarrow \A$ continuously (the union might be discontinuous at $\partial \A_0$), contrary to the smooth case. 
\end{remark}
\begin{proof}[Proof of \cref{Claim E1}] The proof is done basically as for \cref{Claim CB2}. The only change is that the map $\mathring h_n$ defined by \cref{Katok classic} is replaced by the map defined in \cref{defhn} using \cref{neck lace}. 
Then we use \cref{coro:approx analytic} to obtain an analytic map $H_n$ whose restriction to $\A_{\delta_n}$ is close to the one of $H_{n-1}\circ \mathring h_n$. This enables us to obtain statements $(P_1)$ and $(P_2)$ of \cref{Claim E1} from the proof of statements $(P_1)$ and $(P_2)$ in \cref{Claim E0}. 
\end{proof}

\subsection{Proof of the corollaries}\label{proof of corollaries}
\begin{proof}[Proof of \cref{maincoro}]
By \cref{main1}, there is an entire symplectomorphism $F$ of $\A$ which leaves invariant an analytic cylinder $\A_0'\subset \A$.

 1.) For $\pm \in \{-,+\}$, let $\gamma^\pm\in C^\omega(\T,\R)$ be such that $\A'_0=\{(\theta,y)\in \A: \gamma^-(\theta) < y < \gamma^+(\theta)\}.$
Note that up to a conjugacy with a map of the form $(\theta, y)\mapsto (\theta,\alpha \cdot y)$ with $\alpha>0$, we can assume that the mean of $ \gamma^+- \gamma^- $ is 2. 
Now let $\ell(\theta)= (\gamma^+- \gamma^-)(\theta)/2$. 
Note that $\int_\T \ell d\theta= 1$. 

Thus there exists $L\in C^\omega(\T,\T)$ such that $\partial_\theta L= \ell$ and $L(0)=0$. 
Note that $\partial_\theta L$ does not vanish and so $L$ is a diffeomorphism. Set:
\[ 
\psi: (\theta, y)\mapsto (
L (\theta), y/\ell(\theta)).\] 
Note that $\psi$ is analytic and symplectic. Furthermore, there is $\gamma^0\in C^\omega(\T, \R)$ such that the image of 
$\A'_0$ by $\psi$ is $ \{(\theta,y)\in \A: \gamma^0(\theta) -1\le y \le \gamma^0(\theta)+1\}$. Let $h$ be the composition of $\psi$ with the map $(\theta,y)\mapsto (\theta, y-\gamma^0(\theta))$. Observe that $h$ sends $\A_0'$ to $\A_0$. Then $h\circ f\circ h^{-1}$ satisfies the desired properties. 

2) Let $\hat \A_0$ be a neighborhood of $\A_0'$ which is bounded by two invariant curves $\Gamma_1:= \{(\theta, \gamma_1(\theta)): \theta\in \T\}$ and $\Gamma_2:= \{(\theta, \gamma_2(\theta)): \theta\in \T\}$. By \cite[Thm 2.14]{berger2021coexistence}, as the system coincides with the flow map of a non-degenerate Hamiltonian (which is the second coordinate of the conjugacy map) on a neighborhood of $\Gamma_1\sqcup \Gamma_2$, 
we can blow down the latter curves to transform $\hat \A_0$ to a sphere. Then the dynamics descends to an analytic and symplectic diffeomorphism of the sphere.

 Let us give a short independent (sketch of) proof on how to complement $\hat \A_0$ to obtain a sphere. For $i\in \{1,2\}$, let $V_i$ be a small neighborhood of $\Gamma_i$ in which the dynamics is analytically conjugate to the rotation $(\theta,y)\in \R/\Z\times [-\epsilon,\epsilon]\mapsto (\theta+\alpha,y)$ via a map $\phi_i: V_i\to \R/2\pi \Z\times [-\epsilon,\epsilon]$ so that $\phi_i(\Gamma_i)= \R/2\pi \Z\times \{0\}$. Up to a composition with $-id$, we can assume that $\phi_i(V_i\cap \hat \A_0)$ is $\R/2\pi \Z\times [0,\epsilon]$. 
 Now we compose $\phi_i$ with the analytic symplectomorphism $(\theta,r)\in \R/2\pi \Z\times [0,\epsilon]\mapsto (\sqrt{1+2r} \cos\theta,\sqrt{1+2r} \sin\theta)$. From this we obtain that the composition $\psi_i$ is an analytic symplectomorphism from $V_i\cap \hat \A_0$ onto the annulus $C_i$ of radii $(1,\sqrt{1+2\epsilon})$ and which conjugates $f| V_i\cap \hat \A_0$ to a map which coincides with the rotation of angle $\alpha$ on the disk of radius $\sqrt{1+2\epsilon}$. Now we glue the cylinder $\hat \A_0$ to two copies $\D_1$ and $\D_2$ of the disk of radius $\sqrt{1+2\epsilon}$ at $C_1$ and $C_2$ via resp. $\psi_1$ and $\psi_2$. This forms a sphere on which the dynamics extends canonically to the rotation of angle $\alpha$ on the inclusions of $\D_1$ and $\D_2$. This extension is ergodic on the inclusion of $\A'_0$ and has no more than two periodic points. 
 \end{proof} 

\begin{proof}[Proof of \cref{main2}]
Let $F$ be given by Theorems \ref{main1} or \ref{main3}. We recall that $F$ does not have periodic points in $\C/\Z\times \C$ and so that $J^*=\emptyset$. Hence it suffices to show that its Julia set $J$ is nonempty. More precisely, we are going to show that the cylinder $\A_0'$ is included in the Julia set of $F$. 
Let $x\in \A_0'$. Both the forward and backward orbits of $x$ are bounded. 

For the sake of contradiction, assume that there is an open neighborhood $U$ of $x\in \C/\Z\times \C$ such that the following set is bounded:
 \[U^+:= \bigcup_{n\in \N} F^n(U)\; .\]
Note that $F^{-1}(U^+) \supset U^+$. As $F$ preserves the volume, $F^{-1}(U^+)$ and $U^+$ have the same volume. As $F^{-1}(U^+)$ is open, it cannot contain a point 
which is not in the closure of $U^+$. So $F^{-1}(U^+)\subset cl(U^+)$. Likewise, we have $F^{-n}(U^+)\subset cl(U^+)$ for every $n$ and so:
\[\hat U:= \bigcup_{n\ge 0} F^{-n}(U^+)= \bigcup_{n\in \Z} F^n(U)\ \]
is a bounded $F$-invariant invariant open set. 
Then one proves, as in \cite[Appendix]{BSR2} or \cite[Thm 2.1]{Be18}, that $G:= cl(\{F^n|\hat U: n\in \Z\}$ is a compact Abelian subgroup of complex automorphisms of the bounded domain $\hat U$. Then by the Cartan Theorem \cite{Na71}, the set $G$ is a Lie group equal to an $i$-dimensional torus $\T^i$. As the set of iterates of $F|\hat U$ is dense in $G$, the element $F|\hat U$ must act on $G$ as an irrational rotation. Thus the restriction of $F$ to the closure of the orbit of any $x\in \hat U\cap \A$ is semi-conjugate via an analytic map to an irrational rotation on $\T^i$. The rank of the differential of the semi-conjugacy must be invariant by the irrational rotation and so constant. 
Therefore, $G\cdot x=cl \{F^n(x): n\in \Z\}$ must be a torus analytically embedded into $\A$. This torus cannot be of dimension $\ge 2$ (as $\T^2$ cannot be embedded in $\A$) nor of dimension $0$ (as $F$ has no periodic point). Thus the orbit of $x$ is included in an analytic circle (on which $F$ acts as an irrational rotation). 
 This implies that there is an analytic fibration by circles on $\A_0'$ which is left invariant by the dynamics. This is in contradiction to the case of \cref{main1} which asserts that a typical point of $\A_0'$ has the closure of its orbit equal to $cl(\A_0')$. This is also in contradiction to the case of \cref{main3} which asserts that the local order of the emergence is $2$, while by \cref{casestudy}.2, a dynamics leaving invariant such a differentiable fibration should have an order of emergence $0$. \end{proof}

\section{Approximation Theorems}\label{proof approx analytic}
In order to prove \cref{approx analytic} we are going to study the generators of $\Ham^\infty_0(\A)$ and $\Ham^\omega(\A)$, see Theorems \ref{approx reel1} and \ref{approx analytic1}. 
In \cref{sec:proof of approx analytic}, we will prove that \cref{approx analytic1} implies \cref{approx analytic}. 

\subsection{Generators of $\Ham$}
The main theorem of \cite{BT22} implies that $\Ham^\infty (\A)$ is spanned by the following subgroups:
\begin{itemize}
\item Let $\cal V$ be the subgroup of $\Ham^\infty( \A)$ of the form $(\theta, y)\mapsto (\theta , y +v(\theta))$ with $\int vd\theta=0$.
\item 
Let $\cal T$ be the subgroup of $\Ham^\infty( \A)$ of the form $(\theta, y)\mapsto (\theta+\tau (y), y)$. 
\end{itemize}

\begin{theorem}[with Turaev] \label{BT thm}
Any maps of $\Ham^\infty(\A)$ can be arbitrarily well $C^\infty$-approximated by a composition of maps in $\cal T$ or $\cal V$ 
\end{theorem} 
Interestingly, the proof is constructive. Also, as a consequence of the proof using a Lie bracket technique and Fourier's Theorem inspired from \cite{BGH22}, we will show that any maps of $\Ham^\infty(\A)$ can be arbitrarily well approximated by a finite composition of maps in $\cal T$ or $[\cal V,\cal T]$, where:
\[ [\cal V,\cal T]:= \{ [V,T]=V^{-1}\circ T^{-1}\circ V\circ T: V\in \mathcal V\; \&\; T\in \mathcal T \}\; .\]
Also, the subgroup of entire maps of $\cal T$ is $C^\infty$-dense in $\cal T$ and the subgroup of entire maps of $\cal V$ is $C^\infty$-dense in $\cal V$. Yet to prove the main approximation \cref{approx analytic}, we should take care of the supports of the maps in this decomposition. First, we will see that any maps of $\Ham^\infty_0(\A)$ can be arbitrarily well $C^\infty$-approximated by a composition of maps in $\cal T\cap \Ham^\infty_0(\A)$ or $[\cal V , \cal T\cap \Ham^\infty_0(\A)]\cap \Ham^\infty_0(\A)$. Indeed, the set $[\cal V , \cal T\cap \Ham^\infty_0(\A)]$ is not included in $ \Ham^\infty_0(\A) $, so we shall study more carefully the supports of these decompositions by introducing the following notations:
\begin{itemize}
\item For $\epsilon>0$, let $\cal V_\epsilon:= \{V\in {\cal V}: \sup_{x\in \A} \|V(x)-x\|<\epsilon\}$ 
be the $C^0$-$\epsilon$-nghbd of $id$ in $\cal V$.
\item For $\delta>0$, let $\cal T_\delta:= \Ham^\infty_\delta(\A)\cap \cal T $, where $\Ham^\infty_\delta(\A)$ is defined in \cref{def hamdelta}. 
\item Let $ [\cal V_\epsilon,\cal T_\delta]:= \{ [V,T]=V^{-1}\circ T^{-1}\circ V\circ T: V\in \mathcal V_\epsilon\; \&\; T\in \mathcal T_\delta \} $.
\end{itemize}
Most of the maps of $\cal V_\epsilon$ are not compactly supported. In contrast, $\cal T_\delta$ and $ [\cal V_\epsilon,\cal T_\delta]$ are formed by compactly supported maps if $\delta>\epsilon>0$:
\begin{fact} \label{fact reel}For every $\delta>\epsilon>0$, the set $ [\cal V_\epsilon,\cal T_\delta]$ is formed by maps in 
$\Ham^\infty_{\delta-\epsilon }(\A)$. 
\end{fact}
\begin{proof}
Recall that 
$\A\setminus 
\A_{\delta-\epsilon }= \T\times (-\infty, -1+ \delta-\epsilon)\sqcup \T\times(1- \delta+\epsilon,\infty) $. 
Let:
\[ [V,T]=V^{-1}\circ T^{-1}\circ V\circ T\in [\cal V_\epsilon,\cal T_\delta].\] Every $z$ in $\A\setminus 
\A_{\delta-\epsilon }$ is fixed by $T$ and sent by $V$ into $ \A\setminus \A_\delta$ on which $T^{-1}$ is the identity. Thus $T^{-1}\circ V\circ T(z)=V(z)$ and so $[V,T](z)=z$.
\end{proof}
Here is a compactly supported counterpart of \cref{BT thm}, that we will prove in \cref{proof approx reel1}:
\begin{theorem} \label{approx reel1}
For any $0<\epsilon<\delta<1$, any map $F\in \Ham^\infty_\delta(\A)$ and any $C^\infty$-neighborhood $\cal U$ of $F$, there is a composition 
$\tilde F:= F_1\circ \cdots\circ F_M$ of maps $F_j$ in $\cal T_\delta$ or in $[\cal V_\epsilon,\cal T_\delta ]$ such that $\tilde F$ is in $\cal U$. 
\end{theorem} 
It is indeed a compactly supported counterpart of \cref{BT thm} since by \cref{fact reel}, each map $F_j$ belongs to $\Ham_{0}^\infty(\A)$. We shall prove \cref{approx analytic} by introducing the analytic counterpart of \cref{approx reel1}, which requires the following notations. First recall that the sets $K_\rho$ and $\T_\rho$ were defined in \cref{Krho recalled} \cpageref{Krho recalled}.
Recall that $\Ham^\omega (\A)$ denotes the space of entire Hamiltonian maps of $\A$ and $\Ham_{\rho}^\omega (\A)$ is its subset formed by maps whose restrictions to $K_{\rho} $ are $\rho^{-1}$-$C^0$-close to the identity. 

We will use the following generators for $\epsilon,\eta>0$ and $\rho>1$:
\begin{itemize}
\item Let $\cal T_{\rho}^\omega:= \cal T\cap \Ham^\omega_{\rho} (\A)$. In other words, it consists of maps of the form:
\[ (\theta, y)\mapsto (\theta+\tau (y), y)\text{ where $\tau $ is an entire function satisfying }\sup_{Q_{\rho }} |\tau |<\rho^{-1}.\] 
\item Let $\cal V_{\rho}^\omega:= \cal V \cap \Ham_{\rho}^\omega (\A)$. 
In other words, it consists of maps of the form:
\[ (\theta, y)\mapsto (\theta, y+v(\theta))\text{ where $v $ is an entire function satisfying }\sup_{\T_{\rho }} |v |<\rho^{-1}\text{ and } \int_\T v\, d\leb= 0\; .\] 
\item For any $\rho_1,\rho_2>1$, let
$ [\cal V_{\rho_1}\, , \cal T_{\rho_2}]:= \{ [V,T]=V^{-1}\circ T^{-1}\circ V\circ T: V\in \mathcal V_{\rho_1}\; \&\; T\in \mathcal T_{\rho_2} \}\; .$
\end{itemize} 

Here is the analytic counterpart of \cref{approx reel1}, that we will prove in \cref{proof approx analytic1}:

\begin{theorem} \label{approx analytic1}
 For any $0<\delta <1$, for every sequence $(\rho_j)_{j\ge 0}$ of numbers $>1$, any map $F\in \Ham^\infty_\delta(\A)$ and any neighborhood $\cal U$ of $F| \A_\delta $ in $C^\infty(\A_\delta,\A)$, there is a composition $\tilde F:= F_1\circ \cdots \circ F_M$ of maps $F_j$ in $\cal T_{\rho_j}^\omega$ or in $[\cal V_{\rho_0}^\omega ,\cal T_{\rho_j} ^\omega]$ whose restriction $\tilde F| \A_\delta $ is in $\cal U$. 
\end{theorem} 
\begin{remark} 
In \cref{approx analytic1} the approximation is less precise than in \cref{approx reel1}, but it is done by analytic twist maps satisfying furthermore nice bounds.
\end{remark} 
Let us now prove  the approximation \cref{approx reel1,approx analytic1}. 
\subsection{The smooth case: proof of \cref{approx reel1}}\label{proof approx reel1}
For $\delta>0$, let $C^\infty_\delta(\A,\R)$ and $C^\infty_\delta(\R,\R)$
be the spaces of smooth functions with support in respectively $\A_\delta$ and $[\delta-1, 1-\delta]$. We recall that the \emph{symplectic gradient} of $H\in C^\infty(\A,\R)$ is:
\[X_H:= (\partial_y H, -\partial_\theta H)\; .\]

The proof of this theorem will use Poisson brackets. We recall that given two functions $f,g\in C^\infty(\A, \R)$ the \emph{Poisson bracket} $\{f,g\}$ is the function defined by:
\[
\{f,g\}= 
\partial_\theta f \cdot \partial_y g- \partial_\theta g \cdot \partial_y f\; .\]
The Poisson bracket has the following well-known property:
\begin{proposition}\label{rel Lie Poisson bracket} 
Let $H_1, H_2\in C^\infty(\A, \R)$ and $H=\{H_1, H_2\}$ be with respective symplectic gradients denoted by $X_1$, $X_2$ and $X$. Then $X$ equals the Lie bracket of $[X_1,X_2]= DX_2(X_1)-DX_1(X_2)$:
\[ [X_1,X_2] = X\; .\]
 \end{proposition} 

In order to prove \cref{approx reel1} we will use:
\begin{lemma}\label{decomposition fourier}
Let $\delta>0$ and $H\in C^\infty_\delta(\A,\R)$. For every neighborhood $ \cal U$ of $H$ in $C^\infty(\A,\R)$, there exist a 
function $C\in C^\infty_\delta (\R,\R)$ and $2N$ functions $ A_1,\dots, A_N,B_1,\dots, B_N \in C^\infty (\R,\R)$ such that the following is in $\cal U$:
\begin{equation}\label{deftildeH} \tilde H: (\theta, y)\in\A\mapsto C(y) + \sum_{m=1}^N\left\{\frac{ \cos (2\pi m\theta)}{2\pi m} , A_m(y)\right\} +\left \{\frac{ \sin (2\pi m\theta)}{2\pi m} , B_m(y)\right\}\; .\end{equation} 
Moreover, the derivatives of the functions $ A_1,\dots, A_N,B_1,\dots, B_N $ are in $ C^\infty_\delta (\R,\R)$. 
\end{lemma}
\begin{proof}
Using Fourier series decomposition, there are $C\in C^\infty_\delta (\R,\R)$ and $2N$ functions $ a_1,\dots, a_N,b_1,\dots, b_N \in C^\infty_\delta (\R,\R)$ such that the following is in $\cal U$:
\[ \tilde H(\theta, y)\in\A\mapsto C(y) + \sum_{m=1}^N a_m(y) \cos (2\pi m \theta) + b_m(y) \sin (2\pi m \theta) \; .\]
Let $A_m:=- \int_0^y b_m(t) dt$ and $B_m:=\int_0^y a_m(t) dt$ whose derivatives are indeed in $ C^\infty_\delta (\R,\R)$. Note that: 
\[\left\{\frac{ \cos (2\pi m\theta)}{2\pi m} , A_m(y)\right\} = b_m (y)\cdot \sin(2\pi m\theta)\qand \left \{\frac{ \sin (2\pi m\theta)}{2\pi m} , B_m(y)\right\}= a_m (y)\cdot \cos(2\pi m\theta)\; .\]
Thus $\tilde H$ has the desired form. 
\end{proof}
A consequence of the later lemma and \cref{rel Lie Poisson bracket} is: 
\begin{lemma}\label{approx homogenous}
The map $\tilde H$ of Lemma \ref{decomposition fourier} satisfies that there exist functions 
$(\tau_{j})_{ 0\le j\le 2N} \in C^\infty _\delta(\R,\R)$ and $(v_{ j})_{ 1\le j\le 2N} \in C^\infty (\T,\R)$ such that 
with $\phi^t_{\tilde H}$ the time-$t$ map of $\tilde H$ and with:
\[T_{j }^t:=(\theta,y)\mapsto (\theta+t \tau_{j}(y), y) \qand V_{j}^t:=(\theta,y)\mapsto (\theta, y+t v_{j}(\theta)) \] 
 it holds (at $N$ fixed) in any $C^r$-topology when $t\to 0$:
\[ \phi^{t^2}_{\tilde H}= [ T^t_{ 2N}, V^t_{ 2N}] \circ \cdots \circ [ T^t_{ j}, V^t_{ j}]\circ \cdots \circ [ T^t_{ 1}, V^t_{ 1}]\circ T^{t^2}_{ 0}+o(t^2)\; . \]
Moreover, the $C^0$-norms of $v_j$ are bounded by $1$. 
\end{lemma}
\begin{proof}In the setting of \cref{decomposition fourier}, let $H_0:=(\theta, y)\in\A\mapsto C(y) $ and for $1\le m\le N$, put $H_{2m}:=(\theta, y)\in\A\mapsto \{\tfrac1{2\pi m} \cos (2\pi m\theta) , A_m(y)\}$ and $H_{2m-1}:=(\theta, y)\in\A\mapsto \{\tfrac1{2\pi m} \sin (2\pi m\theta) , B_m(y)\}$. 
Let $X_{H_j}$ be the symplectic gradient of $H_j$ and $X_{\tilde H}$ the one of ${\tilde H}$. As ${\tilde H}= \sum_j H_j$, we have:
\[X_{\tilde H}= \sum_j X_{H_j}\; .\]
Hence, denoting $\phi^t_{H_j}$ the flow map of $X_{H_j}$, in any $C^r$-topology $\infty>r\ge 1$, we have when $t\to 0$:
\begin{equation}\label{premiere compo} \phi^{t}_{\tilde H}= id +t X_{\tilde H}+ O(t^2)=id +t \sum_{j=0}^{2N} X_{H_j}+ O(t^2)\; . 
\end{equation} 
Now observe that with $ \tau_0:= C'(y)$ and $T_{0 }^t(\theta,y)=(\theta+t \tau_{0}(y), y)$, it holds $id+t X_{H_0}(\theta,y)= T^t_{0} (\theta,y)
$. Thus
\begin{equation}\label{premiere compobis} \phi^{t}_{\tilde H}= (id +t \sum_{j=1}^{2N} X_{H_j} )\circ T^t_{0} +O(t^2)\; .\end{equation}  
Also, for $1\le j\le 2N$, we have $H_{j}$ of the form $H_j:(\theta, y)\in\A\mapsto \{f_j(\theta) , g_j(y)\}$. Put:
\[v_j:=- f_j'\qand \tau_j:= g_j'\; .\]
Note that $\sup |f'_j| =1$. With 
$T_{j }^t (\theta,y)= (\theta+t \tau_{j}(y), y)$ and $V_{j}^t(\theta,y)= (\theta, y+t v_{j}(\theta))$, the symplectic gradients of $(\theta,y)\mapsto f_j(\theta)$ and $(\theta,y)\mapsto g_j(y)$ are:
\[ X_{f_j} (\theta,y) = (0, v_j(\theta) )=\partial_t V_j(\theta,y) \qand 
X_{g_j} (\theta,y) = (\tau_j(y), 0 )=\partial_t T_j(\theta,y) \; .\] 
Then by \cref{rel Lie Poisson bracket}, the symplectic gradient of $H_j$ is equal to the Lie bracket of the symplectic gradients of $f_j$ and $g_j$:
\[ X_{H_j} = [\partial_t T^t_j, \partial_t V^t_j] \; .\] 
Thus we have in any $C^r$-topology:
\[ id+ t^2 X_{H_j} = [ T^t_j, V^t_j] +O(t^3)\; .\] 
Injecting the latter into \cref{premiere compobis} at time $t^2$, it comes:
\begin{equation}\label{premiere compoters} \phi^{t^2}_{\tilde H}= 
[ T^t_{2N}, V^t_{2N}]\circ \cdots \circ [ T^t_{1}, V^t_{1}]\circ T^{t^2}_{0} +O(t^3)\end{equation} 
which is the desired result. 
\end{proof}

\begin{proof}[Proof of \cref{approx reel1}]
For any $f\in \Ham^\infty_\delta (\A)$ , there is a smooth family $(H_t)_{t\in [0,1]}$ of Hamiltonians $H_t\in C^\infty_\delta (\A)$ which defines a family $(f_t)_{t\in [0,1]}$ such that $f_0=id$, $f_1=f$ and $\partial_t f_t$ is the symplectic gradient of $H_t$. Observe that for $M$ large, we have:
\[f= (f\circ f_{(M-1)/M}^{-1})\circ \cdots \circ (f_{i/M}\circ f_{(i-1)/M}^{-1})\circ 
\cdots \circ f_{1/M}\; .\]
Each $(f_{j/M}\circ f_{(j-1)/M}^{-1})$ is $O(M^{-2})$-$C^r$-close to the time $\tau=M^{-1}$ map $\phi_j^{\tau}$ of the vector field $\partial_t f_{j/M}$. Thus in the $C^r$-topology:
\[f= \phi_M^{\tau }\circ \cdots \circ \phi_1^{\tau}+O(M^{-1})\; .\]
Note that $\phi_j^{\tau }$ is the time $\tau $ map of the symplectic gradient of $H_{j/M}\in C^\infty_\delta (\A)$. 
By Lemmas \ref{decomposition fourier} and \ref{approx homogenous}, with $t^2:= \tau$, each $\phi_j^{\tau }$ is equal to a composition of elements in $\cal T_\delta$ and $[\cal V_{t} ,\cal T_\delta]$ up to precision $o(\tau)=o(M^{-1})$. 
Thus $f$ is $C^r$-close to a composition of elements in $\cal T_\delta$ and $[\cal V_{t},\cal T_\delta]$ for any $t>0$. We conclude by choosing $t^2<\epsilon$. \end{proof}
\subsection{The analytic case: proof of \cref{approx analytic1}}\label{proof approx analytic1}
The proof of \cref{approx analytic1} follows basically the same lines as the one of \cref{approx reel1}, although we have to extend the bounds to a complex compact set and use the Runge theorem.

In this subsection, we fix $1>\delta >0$ and a sequence $(\rho_j)_{j\ge 0}$ of numbers $>1$ as in the statement of \cref{approx analytic1}. For every $n\ge 1$, put: 
\[M_n := \rho_0\cdot \sup_{z\in \T_{\rho_0}} |2\pi n\cdot \exp(2\pi i\cdot n \cdot z)|\; . \]
Then observe that the same proof as \cref{decomposition fourier} gives:
 \begin{lemma}\label{decomposition fourier2}
 The same statement as \cref{decomposition fourier} is true with the following formula instead of \cref{deftildeH}:
 \[ \tilde H: (\theta, y)\in\A\mapsto C(y) + \sum_{m=1}^N\left\{\frac{ \cos (2\pi m\theta)}{M_m} , A_m(y)\right\} +\left \{\frac{ \sin (2\pi m\theta)}{M_m} , B_m(y)\right\}\; .\]
 \end{lemma}
Then by plugging \cref{decomposition fourier2} instead of \cref{decomposition fourier} in the proof of \cref{approx homogenous} we obtain:
\begin{lemma}\label{approx homogenous ana}
The map $\tilde H$ of Lemma \ref{decomposition fourier2} satisfies the same property as in \cref{approx homogenous} and moreover 
 each $v_j$ can be chosen entire and satisfying $\sup\{ |v_j(z)|:{z\in \T_{\rho_0}}\} \le 1/\rho_0$. 
\end{lemma}
Using the latter lemma instead of \cref{approx homogenous} in the proof of \cref{approx reel1} gives immediately:
\begin{lemma} \label{approx reel1 ana}
Any $F\in \Ham_\delta^\infty(\A)$ can be arbitrarily well $C^\infty$-approximated by a composition 
$F_M\circ \cdots \circ F_1$ with each $F_j$ in $\cal T_\delta$ or in $[\cal V_{\rho_0}^\omega,\cal T_\delta ]$. 
\end{lemma} 
We are now ready for the:
\begin{proof}[Proof of \cref{approx analytic1}]
Let $1>\delta >0$ and let $(\rho_j)_{j\ge 0}$ be a sequence of numbers $>1$ as in the statement of the theorem. As 
$\cal V^\omega_\rho\subset \cal V_{\rho'}^\omega$ when $\rho'\le \rho$, we can assume:
\begin{equation}\label{cond delta epsilon2}2\rho^{-1}_0<\delta \; .\end{equation}

Let $F\in \Ham^\infty_\delta(\A)$ and let $\cal U$ be a neighborhood of the restriction of $F|\A_\delta$ in the $C^\infty$-topology.

 By \cref{approx reel1 ana}, there exists a composition $F_M\circ \cdots \circ F_1 $ whose restriction 
 to $\A_\delta$ is in $\cal U$ and formed by maps $F_j$ in $\cal T_\delta$ or in $[\cal V_{\rho_0}^\omega,\cal T_\delta ]$. It suffices to change the maps $F_j$ in $\cal T_\delta$ to maps in $\cal T_{\rho_j}^\omega$ and those in $[\cal V_{\rho_0}^\omega,\cal T_\delta ]$ to maps in $ [\cal V_{\rho_0}^\omega ,\cal T_{\rho_j}^\omega ]$ so that the restriction of their compositions to $\A_\delta= \T\times [-1+ \delta,1-\delta] $ remains in $\cal U$. 
 
 We recall that the support of each $F_j$ is in $\A_{\delta-\rho_0^{-1}}$ by \cref{fact reel}.
Thus $F_j$ sends any small neighborhood of $\A_{\delta-\rho_0^{-1}}$ into another small neighborhood of $\A_{\delta-\rho_0^{-1}}$. Likewise, if $\tilde F_j$ is a map whose restriction to a neighborhood of $\A_{\delta-\rho_0^{-1}}$ is close to the one of $F_j$, then $\tilde F_j$ sends any small neighborhood of $\A_{\delta-\rho_0^{-1}}$ to another small neighborhood of $\A_{\delta-\rho_0^{-1}}$.
 Thus the composition $\tilde F_M\circ \cdots \circ \tilde F_1$ has its restriction to $\A_\delta \subset \A_{\delta-\rho_0^{-1}}$ close to the one of $F_M\circ \cdots \circ F_1 $ and so close to $F|\A_\delta$. 
Hence it suffices to show:
\begin{claim}\label{dernier claim a adapter} For every $j $, there exists $\tilde F_j$ in $ \cal T_{\rho_j}^\omega\cup [\cal V_{\rho_0}^\omega ,\cal T_{\rho_j}^\omega ]$ whose restriction to a neighborhood of $\A_{\delta-\rho_0^{-1}}$ is $C^\infty$-close to the one of $F_j$. \end{claim}

To show this claim, we fix a neighborhood $\cal N$ of $[-1+\delta-2\rho_0^{-1},1-\delta +2\rho_0^{-1}]$ which is disjoint from $\R\setminus (-1,1)$. Such a neighborhood exists by \cref{cond delta epsilon2}. There are two cases.
 
\underline{Case 1)} If $F_j\in \cal T_{\delta}$, then there exists $\tau_j$ such that $F_j:(\theta,y)\mapsto (\theta+\tau_j(y), y)$. First, by the Weierstrass Theorem, there exists a polynomial $C^\infty$-approximation $\hat \tau_j$ of $\tau_j|\cal N $. Now let $\tilde {\cal N}$ be a neighborhood of $cl(\cal N)$ in $\C\setminus Q_{\rho_j}$ so that the complement of $\tilde {\cal N}\cup Q_{\rho_j}$ is connected. See \cref{figure_Runge}.
\begin{figure}[h]
\centering
\includegraphics[width=7cm]{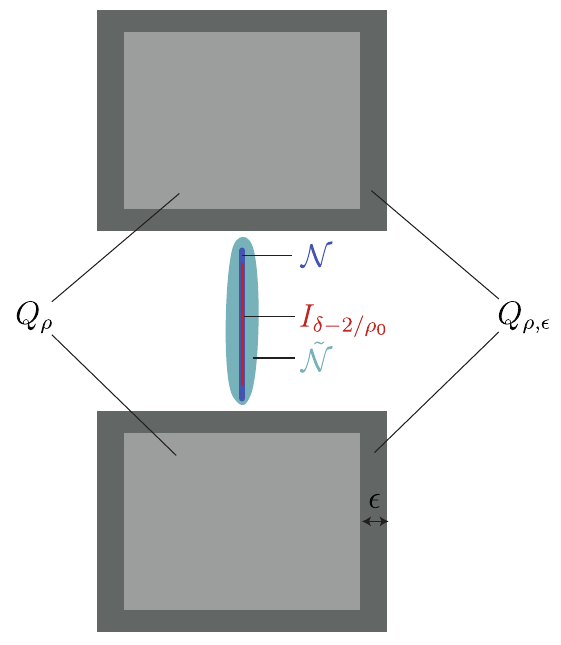}
\caption{Domains of $\C$ involved in the application of the Runge theorem.}\label{figure_Runge}
\end{figure}

 Then, using Runge's theorem, there exists an entire map $\tilde \tau_j$ which is $C^0$-close to $\hat \tau_j$ on $\tilde {\cal N} $ and which is $\rho^{-1}_j$-small on $Q_{\rho_j}$. As $\hat \tau_j$ is real (i.e. it commutes with $z\mapsto \bar z$), up to considering $\tfrac12 \big(\tilde \tau_j (z)+ \overline {\tilde \tau_j (\bar z)}\big)$ instead of $\tilde \tau_j$, we can assume that $\tilde \tau_j$ is real. Then by the Cauchy inequality, the map $\tilde \tau_j$ restricted to $\cal N$ is $C^\infty$-close to $\hat \tau_j|\cal N$ and so $ \tau_j|\cal N$. 
 
 Then we observe that $\tilde F_j:(\theta,y)\mapsto (\theta+\tilde \tau_j(y), y)$ belongs to $\cal T_{\rho_j}^\omega$ and has its restriction to $\T\times \cal N $ close to the one of $F_j$.
 
\underline{Case 2)} If $F_j\in [\cal V_{\rho_0}^\omega,\cal T_\delta ]$, then there exist $V_j\in \cal V_{\rho_0}^\omega$ and $T_j:(\theta,y)\mapsto (\theta+\tau_j(y), y)$ such that $F_j=[V_j,T_j]$. As in the case $(1)$, there is $\tilde T_j\in \cal T_{\rho_j}^\omega$ such that the restriction of $T_j$ and 
 $\tilde T_j$ to $\T\times \cal N $ are $C^\infty$-close. Then put:
\[\tilde F_j= [V_j,\tilde T_j] = V_j^{-1}\circ \tilde T_j^{-1} \circ V_j\circ \tilde T_j \; .\] 
Note that a small neighborhood of $\A_{\delta-\rho_0^{-1}}$ is sent by $\tilde T_j$ to a small neighborhood of $\A_{\delta-\rho_0^{-1}}$,
 which is then sent by $V_j$ into a small neighborhood of $\A_{\delta-2\rho_0^{-1}}$, which is therefore included in $\T\times \cal N $. 
 As $\tilde T_j^{-1}|\T\times \cal N $ is close to $T_j^{-1}|\T\times \cal N $, it comes that the restrictions of $\tilde F_j$ and $F_j$ to a neighborhood of $\A_{\delta-\rho_0^{-1}}$ are $C^\infty$-close. 
\end{proof}

\subsection{ Proof that \cref{approx analytic1} implies \cref{approx analytic} }\label{sec:proof of approx analytic} 
By \cref{approx analytic1}, we know that every map $F\in \Ham^\infty_\delta(\A)$ has its restriction to $\A_\delta$ which can be approximated by a composition $\tilde F:= F_1\circ \cdots \circ F_M$ of maps $F_j$ in $\cal T^\omega_{\rho_j}\cup [\cal V^\omega_{\rho_0} ,\cal T^\omega_{\rho_j} ]$, for any chosen sequence $(\rho_j)_j$ of numbers $>1$. Yet we do not know if 
$\tilde F|K_\rho$ is small. Actually, it might not be  the case since the set $K_\rho$ is a priori not left invariant by the maps $F_j$. Thus we shall first deduce from \cref{approx analytic1} a similar statement (see \cref{approx analytic2}) giving an approximation of $F|\A_\delta$ by the restriction of a composition of maps $F_j$ which are moreover small on a neighborhood $K_{\rho,\epsilon}$ of $K_\rho$. Then we will show that the composition of these maps is indeed small on $K_\rho$.

This leads us to introduce a few technical notations. 

\paragraph{Notations} 
For $\epsilon>0$, we denote $K_{\rho,\epsilon}$ and $Q_{\rho,\epsilon}$ the $\epsilon$-neighborhoods of $K_\rho $ and $Q_\rho$:
\[ K_{\rho,\epsilon}= \T_{\rho+\epsilon}\times Q_{\rho,\epsilon} \qand Q_{\rho,\epsilon}:= [-\rho-\epsilon,-1+\epsilon]\sqcup [1-\epsilon,\rho+\epsilon] + i [-\rho-\epsilon ,\rho+\epsilon]\; .\]
Similarly, we denote $\Ham_{\rho, \epsilon }^\omega (\A)$ the subset of $\Ham_{\rho}^\omega (\A)$ formed by maps whose restrictions to $K_{\rho,\epsilon} $ is $\rho^{-1}$-$C^0$-close to the identity:
 \[\Ham_{\rho, \epsilon }^\omega (\A):=\left\{F\in \Ham_{\rho }^\omega (\A): \; 
 \sup_{x\in K_{\rho,\epsilon} } |F(x )-x|<\rho^{-1}\right\}\;.\]
Likewise we denote:
\begin{itemize}
\item Let $\cal T_{\rho, \epsilon}^\omega:= \cal T\cap \Ham_{\rho,\epsilon} ^\omega(\A)$. In other words, it consists of maps of the form:
\[ (\theta, y)\mapsto (\theta+\tau (y), y)\text{ where $\tau $ is an entire function satisfying }\sup_{Q_{\rho,\epsilon }} |\tau |<\rho^{-1}.\] 
\item Let $\cal V_{\rho,\epsilon}^\omega:= \cal V \cap \Ham_{\rho,\epsilon}^\omega (\A)$. 
In particular, its elements are of the form:
\[ (\theta, y)\mapsto (\theta, y+v(\theta))\text{ where $v $ is an entire function satisfying }\sup_{\T_{\rho+\epsilon }} |v |<\rho^{-1}.\] 
\item For any $\rho_1,\rho_2>1$, let
$
[\cal V_{\rho_1,\epsilon}^\omega\, ;\cal T_{\rho_2,\epsilon}^\omega]:= \{ [V,T]=V^{-1}\circ T^{-1}\circ V\circ T: V\in \mathcal V_{\rho_1,\epsilon}^\omega\; \&\; T\in \mathcal T_{\rho_2,\epsilon}^\omega \}$.
\end{itemize} 
Similarly to \cref{approx analytic1} we have:
\begin{lemma} \label{approx analytic2}
For any $1>\delta >\epsilon>0$, for every sequence $(\rho_j)_{j\ge 0}$ of numbers $>1$, 
 for every map $F\in \Ham^\infty_\delta(\A)$ and any $C^\infty$-neighborhood $\cal U$ of $F| \A_\delta $, there is a composition $\tilde F:= F_1\circ \cdots \circ F_M$ of maps $F_j$ in $\cal T^\omega_{\rho_j,\epsilon}$ or in $ [\cal V^\omega_{\rho_0,\epsilon} ,\cal T^\omega_{\rho_j,\epsilon} ]$ such that $\tilde F|\A_\delta \in \cal U$. 
\end{lemma} 
\begin{proof} We can assume $\rho_0$ sufficiently large so that: 
\begin{equation}\label{cond delta epsilon3} 2\rho_0^{-1}+\epsilon< \delta\; .\end{equation}
We notice that $\cal V_{\rho_0+\epsilon}^\omega=\cal V_{\rho_0,\epsilon}^\omega$. Thus 
by \cref{approx reel1 ana} with $\rho_0+\epsilon$ instead of $\rho_0$, any $F\in \Ham_\delta^\infty(\A)$ can be arbitrarily well $C^\infty$-approximated by a composition $F_M\circ \cdots \circ F_1$ with each $F_j$ in $\cal T_\delta$ or in $[\cal V_{\rho_0,\epsilon}^\omega,\cal T_\delta ]$. Likewise, it suffices to show:
\begin{claim}\label{dernier claim adapte} For every $j $, there exists $\tilde F_j$ in $ \cal T_{\rho_j,\epsilon}^\omega\cup [\cal V_{\rho_0,\epsilon}^\omega ,\cal T_{\rho_j,\epsilon}^\omega ]$ whose restriction to a neighborhood of $\A_{\delta-\rho_0^{-1}}$ is $C^\infty$-close to the one of $F_j$. \end{claim}
 To prove this claim, we proceed as for \cref{dernier claim a adapter}. We chose a neighborhood $\cal N$ of $[-1+\delta-2\rho_0^{-1},1-\delta +2\rho_0^{-1}]$ which is disjoint from $\R\setminus (-1+\epsilon,1-\epsilon)$. Such a neighborhood exists by \cref{cond delta epsilon3}. 
Then we distinguish two cases. 

\underline{Case 1)} If $F_j\in \cal T_{\delta}$, then it is of the form $F_j:(\theta,y)\mapsto (\theta+\tau_j(y), y)$. First, by the Weierstrass Theorem, there exists a polynomial $C^\infty$-approximation $\hat \tau_j$ of $\tau_j|\cal N $. Now let $\tilde {\cal N}$ be a neighborhood of $cl(\cal N)$ in $\C\setminus Q_{\rho_j,\epsilon}$ so that the complement of $\tilde {\cal N}\cup Q_{\rho_j,\epsilon}$ is connected. See \cref{figure_Runge}. Then, using Runge's Theorem, there exists a real entire map $\tilde \tau_j$ which is $C^0$-close to $\hat \tau_j$ on $\tilde {\cal N} $ and which is $\rho^{-1}_j$-small on $Q_{\rho_j}$. By the Cauchy inequality, the map $\tilde \tau_j$ restricted to $\cal N$ is $C^\infty$-close to $\tau_j$. 

\underline{Case 2)} If $F_j\in [\cal V_{\rho_0,\epsilon}^\omega,\cal T_\delta ]$, then there exist $V_j\in \cal V_{\rho_0,\epsilon}^\omega$ and $T_j:(\theta,y)\mapsto (\theta+\tau_j(y), y)$ such that $F_j=[V_j,T_j]$. As in the case $(1)$, there is $\tilde T_j\in \cal T_{\rho_j,\epsilon}^\omega$ such that the restrictions of $T_j$ and 
 $\tilde T_j$ to a neighborhood of $\T\times \cal N $ are $C^\infty$-close. Then we conclude exactly as in \cref{dernier claim a adapter}. 
\end{proof}   
 A second ingredient for the proof of \cref{approx analytic} is the following counterpart of \cref{fact reel}:
\begin{lemma}\label{fact analytic}
For every $\epsilon>\epsilon'>0$ and $ \rho_k>\rho>1$ such that 
$ \epsilon>\epsilon'+2 \rho^{-1} >0$, for all $[V,T]\in [\cal V^\omega_{\rho, \epsilon}\, ;\, \cal T^\omega_{\rho_k,\epsilon}] $, 
it holds: 
\[\sup_{K_{\rho,\epsilon'}} \|[V,T]-id\|\le \tfrac1 {\rho_k'}\quad \text{with }\rho_k':= \frac{\epsilon-\epsilon'-2\rho^{-1}}{\epsilon-\epsilon'-\rho^{-1}} \cdot \rho_k\, \; .\] 
\end{lemma}
\begin{proof}  
 Let $[V,T]\in [\cal V^\omega_{\rho,\epsilon},\cal T^\omega_{\rho_k,\epsilon}] $. As $\rho_k>\rho$ and $\epsilon>\epsilon'+2\rho^{-1}$, we have:
\begin{equation}\label{recall ineq} \T_{\rho+\epsilon'+2\rho_k^{-1}}\times Q_{\rho, \epsilon'+2\rho^{-1}}\subset K_{\rho,\epsilon'+2\rho^{-1}}\subset K_{\rho,\epsilon}\; .\end{equation}
Note that $T$ sends $ K_{\rho,\epsilon'} =\T_{\rho+\epsilon'}\times Q_{\rho, \epsilon'}$ into $ \T_{\rho+\epsilon'+\rho_k^{-1}}\times Q_{\rho, \epsilon'}$ which is then sent by $V$ into $ \T_{\rho+\epsilon'+\rho_k^{-1}}\times Q_{\rho, \epsilon'+\rho^{-1}}$ and finally sent into $ \T_{\rho+\epsilon'+2\rho_k^{-1}}\times Q_{\rho, \epsilon'+\rho^{-1}}$ by $T^{-1}$.  
Thus, \cref{recall ineq}, the following restrictions are equal:
\[ [V,T]_{|K_{\rho,\epsilon'}}= 
(V ^{-1})_{| K_{\rho,\epsilon'+2\rho^{-1}}}
\circ 
(T^{-1})_{| K_{\rho,\epsilon'+2\rho^{-1}}} 
\circ 
V_{ |K_{\rho,\epsilon'+2\rho^{-1}}}
\circ 
T_{|K_{\rho,\epsilon'}}\; .
\]
As $V$ is a twist map, the restrictions to $K_{\rho,\epsilon}$ of $V$ and its inverse are at equal $C^0$-distance to the canonical inclusion $K_{\rho,\epsilon}\hookrightarrow \A_\C $. This distance is at most $\rho^{-1}$ as $V\in \cal V^\omega_{\rho,\epsilon}$. 
 Hence, by the Cauchy inequality, the restriction of $V^{-1}$ to $K_{\rho,\epsilon'+2\rho^{-1}}$ is 
$C^1$-$\rho^{-1}/ (\epsilon-\epsilon'-2\rho^{-1})$-close to the identity. Hence it holds:
\[\sup_{ K_{\rho,\epsilon'}} \|V^{-1}
\circ 
 T^{-1} 
\circ 
V \circ T
-T\|= \sup_{ T(K_{\rho,\epsilon'})} \|V^{-1}
\circ 
 T^{-1} 
\circ 
V 
-id\| =
\sup_{ V\circ T(K_{\rho,\epsilon'})} \|V^{-1}
\circ 
 T^{-1} 
-V^{-1}\|\]
 \[= 
\sup_{ V\circ T(K_{\rho,\epsilon'})}
\|(V^{-1})_{ |K_{\rho,\epsilon'+2\rho^{-1}}}
\circ 
 T^{-1}
-V^{-1}_{ |K_{\rho,\epsilon'+2\rho^{-1}}}\|
\le \sup \| D(V^{-1})_{| K_{\rho,\epsilon'+2\rho^{-1}}}\| \cdot 
\sup_{ V\circ T(K_{\rho,\epsilon'})}\|T^{-1}-id\|\]
\[
\le \frac{\rho^{-1}}{\epsilon-\epsilon'-2\rho^{-1}}\cdot \sup_{ V\circ T(K_{\rho,\epsilon'})}\|T^{-1}-id\|
\le \frac{\rho^{-1}}{\epsilon-\epsilon'-2\rho^{-1}}\cdot \rho_k^{-1}\; .\]
 Now we infer that $\sup_{K_{\rho_k,\epsilon}}\|T-id\|<\rho_k^{-1}$ to obtain:
 \[
 \sup_{ K_{\rho,\epsilon'}} \|[V,T]-id\| 
 \le 
 \sup_{ T( K_{\rho,\epsilon'})} \|
 V^{-1}
\circ 
T^{-1}
\circ 
V-id\|+
 \sup_{ K_{\rho,\epsilon'}} \|T-id\| \le \frac{\rho^{-1}}{\epsilon-\epsilon'-2\rho^{-1}}\cdot \rho_k^{-1}+ \rho_k^{-1}
\]
which is $\rho_k'^{-1}$. 
\end{proof}
 
 Here is the last lemma needed for the proof of \cref{approx analytic}:
\begin{lemma} \label{fact compo}
For any $\rho >1$ and $\epsilon'> 0$, let $\rho'_k:=2^k\max \{ \epsilon'^{-1} , \rho\} $ for every $k\ge 1$. Then for any sequence $( F_k)_{k\ge 1}$ of maps $ F_k \in \Ham^\omega (\A)$ such that $\sup_{K_{\rho,\epsilon'}} \| F_k-id\|<1/\rho_k'$, it holds:
\[ F_1\circ \cdots \circ F_n\in \Ham_{\rho }^\omega (\A)\text{ for every }n\ge 1.\]
\end{lemma} 
\begin{proof}
Let us first prove the case $n=2$. 
First note that $ F_2(K_{\rho })$ is in $ K_{\rho, 1/\rho_2'}\subset K_{\rho,\epsilon'/4 }\subset K_{\rho,\epsilon'}$. Thus:
\[
\sup_{K_{\rho}}\| F_1\circ F_2-id\| \le 
\sup_{ F_2(K_{\rho})}\| F_1 - id\| +
\sup_{K_{\rho }}\| F_2-id\| \le \sup_{ K_{\rho,\epsilon'} }\| F_1 - id\| +
\sup_{K_{\rho}}\| F_2-id\| 
\]
\[\le 1/\rho'_1+1/\rho'_2\le \rho^{-1} \cdot (2^{-1} +2^{-2})<\rho^{-1}\; .\]
This proves that $ F_1\circ F_2$ belongs to $\Ham_{\rho }^\omega (\A)$. 

Now let us prove the general case for any $n\ge 2$. First note that $ F_{i}(K_{\rho, \sum_{k>i}2^{-k}\epsilon' })$ is included in 
$K_{\rho, \sum_{ k\ge i}2^{-k}\epsilon' }$. Thus, by induction, $ F_i\circ \cdots \circ F_n (K_{\rho})$ is included in 
$K_{\rho, \sum_{ k\ge i}2^{-k}\epsilon' } \subset K_{\rho, \epsilon' } $. We have:
\[
\sup_{K_{\rho}}\| F_1\circ \cdots \circ F_n-id\| \le 
\sum_{0\le i\le n-1} \sup_{ F_i\circ \cdots \circ F_n (K_{\rho})}\| F_{i+1} -id\| \le
\sum_{0\le i\le n-1} \sup_{ K_{\rho,\epsilon'}}\| F_{i+1} -id\|\]
\[\le \sum_{0\le i\le n-1} \rho'^{-1}_{i+1}\le 
\sum_{0\le i\le n-1} 2^{-i-1} \cdot \rho^{-1} <\rho^{-1} \; .
\] 
This proves that $ F_1\circ \cdots \circ F_n$ belongs to $\Ham_{\rho }^\omega (\A)$. 
\end{proof} 
 
 \begin{proof}[Proof of \cref{approx analytic}]
 Let $0<\delta < 1 $ and $ \rho>1$. Let $F\in \Ham_{\delta}^\infty(\A)$ and $\cal U$ a $C^\infty$-neighborhood of the restriction $F| \A_{\delta} $. Let us show that there exists $\tilde F\in \Ham^\omega_{\rho} (\A)$ such that $\tilde F | \A_{\delta} $ is in $\cal U$. 
 
We can assume $\rho>2/\delta$ and so that there are $\delta>\epsilon>\epsilon'>0 $ satisfying $ \epsilon>\epsilon'+2 \rho^{-1}$. For every $k\ge 1$, let:
 \[\rho_k':= 2^k\max \{ \epsilon'^{-1} , \rho\}\qand \rho_k= \frac{\epsilon-\epsilon'-\rho^{-1}}{\epsilon-\epsilon'-2\rho^{-1}}\cdot \rho_k'\; .\]
 
 By \cref{approx analytic2}, there is a composition $\tilde F:= F_1\circ \cdots \circ F_M$ of maps $F_j$ in $\cal T_{\rho_j,\epsilon}^\omega$ or in $ [\cal V_{\rho,\epsilon}^\omega ,\cal T_{\rho_j,\epsilon}^\omega ]$ such that $\tilde F|\A_\delta \in \cal U$. Thus by \cref{fact analytic}, each map $F_j$ satisfies:
 \[\sup_{K_{\rho,\epsilon'}}\|F_j-id\|< \rho_k'^{-1}\; .\]
 Then by \cref{fact compo}, we conclude that the composition $\tilde F$ is in $\Ham_{\rho}^\omega (\A)$. \end{proof} 

\section{Smooth Lemmas}
\subsection{A consequence of Moser's trick}
Let $(S,\omega)$ be a compact surface with boundary. We recall that $ \Sym^\infty_0 (S)$ denotes the space of 
 smooth symplectomorphisms  supported by $S\setminus \partial S$. The following folklore theorem will be a key ingredient for the proof of \cref{Katok classic,neck lace}:
\begin{theorem}\label{coro moser} 
Let $(D_i)_{1\le i\le N} $ and $(D'_i)_{1\le i\le N} $ be two families of disjoint smooth\footnote{i.e. diffeomorphic to a closed disk.} closed disks in $int\, S$, such that for every $i$, $\leb(D_i)=\leb(D_i')$. Then there exists $f\in \Sym^\infty_0 (S)$ such that $f(D_i)=D_i'$ for every $1\le i\le N$. 
\end{theorem} 
We were not able to find a proper reference for a proof of this lemma: Katok \cite[Basic Lemma \textsection 3]{Ka73} and Anosov-Katok \cite[thm 1.3]{AK70} show versions of this theorem which are weaker in dimension 2. Yet they wrote that A. B. Krygin had a proof of such a result (without reference). Others would identify this theorem as a direct application of Moser's trick. \cref{coro moser} is proved by induction using: 
\begin{lemma}\label{moser:lem}
For any subsets $D,D'\subset int\, S$, if $D$ and $D'$ are smooth closed disks and $\leb(D)=\leb(D')$, then there exists $f\in \Sym^\infty_0 (S)$ such that $f(D)=D'$. 
\end{lemma}
\begin{proof}
First we consider an isotopy $(g_t)_{t\in [0,1]}$ of $S$ with support in the interior of $S$ and such that $g_1(D)=D'$, this isotopy can be constructed by retracting the disks to tiny ones and moving one to the other. Each $g_t$ is a diffeomorphism that we shall transform to a symplectic map. Put $D_t:= g_t(D)$.  Using a tubular neighborhood of $\partial D_t$, it is easy to deform smoothly $(g_t)_{t\in [0,1]}$ so that $g_t$ preserves the volume form nearby $\partial D_t$. Indeed, it suffices to work in a tubular neighborhood of $D_t$ and tune the size of the normal component of $g_t$. 

Now we would like to construct a diffeomorphism $\psi_1$ of $S$  which coincides with the identity on a neighborhood of $\partial D\sqcup \partial S$ and such that $g_1\circ \psi_1$ is a symplectomorphism. In order to do so, it suffices to define the restriction of $\psi_1$ to  the closure $N$ of each component of $S\setminus \partial D$.  Let $N_t:= g_t(N)$.  Again, using a tubular neighborhood of $\partial N_t$, it is easy to deform   $(g_t|N)_{t\in [0,1]}$ to a smooth family $(\tilde g_t)_t$ such that 
$\tilde g_t(N)=N_t$, $\tilde g_o= id_N$,  $\tilde g_1= g_1|N$ and  on the neighborhood of $\partial N$ in $N$,  the form $\tilde g_t^* \omega$  coincides with $ \frac{\leb \, N_t} {\leb \,  N}\cdot \omega$.  We define the following form on $N$:
\[\omega_t:=   \frac{\leb \,  N} {\leb  N_t}  \tilde g_t^* \omega .\]
 We notice that $\int_N  \omega_t= \leb(N)$ for every  $t$.  Also $\omega_t$ coincides with $\omega$ on a neighborhood of $\partial N$. Then by Moser's trick \cite[Ex. 3.2.6]{mcduff2017introduction}, there exists a smooth isotopy $\psi_t: N\to N$ such that: $\psi_0=id$ and $\psi_t^*\omega_t = \omega$ and $\psi_t$ is equal to the identity on a neighborhood of $\partial N$. Then 
\[(\tilde g_t\circ \psi_t)^*\omega= \psi_{t}^*( \tilde g_{t}^* \omega)= \frac{\leb \,  N_t} {\leb  N} \psi_{t}^*\omega_t= \frac{\leb \,  N_t} {\leb  N}  \omega\; .\]
This implies that   $(g_1\circ \psi_1)^*\omega=\omega$ on $N$.  
%
%
%
\end{proof}
\begin{proof}[Proof of \cref{coro moser}]
By induction on $n\ge 1$, we construct a map $f_n\in \Sym^\infty_0 (S)$ which sends $(D_i)_{1\le i\le n} $ to $(D'_i)_{1\le i\le n} $. The step $n=1$ is \cref{moser:lem}. Let $n>1$. Let $S_n:= S\setminus \bigcup_{i<n} D'_i$. By \cref{moser:lem}, there exists a map $f\in \Sym^\infty_0 (S_n)$ which sends $f_{n-1}(D_n)$ to $D'_n$. We observe that $f_{n}:= f\circ f_{n-1}$ satisfies the induction hypothesis. 
\end{proof}
 \subsection{Lemma for Anosov-Katok's theorem}\label{proof Katok classic}
\cref{Katok classic} states that: 
 \begin{lemma}For every $q\ge 1$ and $\epsilon>0$, there exists a map $h\in \Ham^\infty_0(\A)$ such that:\begin{enumerate}
\item the map $h$ is $(1/q,0)$-periodic,
\item for every $(\theta,y)\in \A_\epsilon$ the map $ h$ sends the one-dimensional measure $\Leb_{\T \times \{y\}}$ to a measure which is $\epsilon$-close to $\Leb_{\A_0}$. 
\end{enumerate}
 \end{lemma}
 \begin{proof}As both $\Leb_{\A_0}$ and each $\Leb_{\T \times \{y\}}$ are $(1/q,0)$-periodic for every $q\ge 1$, it suffices to show the lemma in the case $q=1$. Indeed, then any $q$-covering of $h$ will satisfy the lemma.
 
Hence for $\epsilon>0$, it suffices to construct $h\in \Ham^\infty_0(\A)$ such that for every $y\in \I_{\epsilon}:= [-1+\epsilon,1-\epsilon]$, the map $ h$ sends the one-dimensional Lebesgue probability measure $\Leb_{\T \times \{y\}}$ to a measure which is $\epsilon$-close to $\Leb_{\A_0}$. The construction is depicted in \cref{figure_ergodic}. 

For $N=n^2$ large, we consider the following collection of boxes $(C_k)_{0\le k \le N-1}$ :
\[C_k:= [\tfrac {k+\epsilon/8}N, \tfrac {k+1-\epsilon/8} N ]\times \I_{\epsilon}\; .\]
 Note that the family $(C_k)_k$ is disjoint  and that for every $k$ we have:
\[\leb_{\A_0} C_k= \frac 2 N(1-\epsilon)(1-\epsilon/4)/\leb\A_0= \frac1N (1-\epsilon)(1-\epsilon/4)\; .\] 
Thus for every $v>\tfrac1N (1-\epsilon)(1-\epsilon/4)$ close  enough to $\tfrac1N (1-\epsilon)(1-\epsilon/4)$, there exists a disjoint family $(D_k)_{1\le k\le N}$ such that each $D_k$ is a neighborhood $C_k$ in $\A_0$,  which is diffeomorphic to a closed disk and whose volume is 
 $\leb_{\A_0} C_k= v$. 
 
 Now recall that $N=n^2$. So there is a bijection $k\in \{0,\dots, N-1\}\to (x_k,y_k)\in \{0,\dots, n-1\}^2$. Let $z_k:= (x_k/n\mod 1, -1+2\cdot y_k/n)\in \A_0$. Note that for $N$ large enough, the measure $N^{-1} \sum_{0\le k<N} \delta_{z_k}$ is $\tfrac \epsilon 2$-close to $\Leb_{\A_0}$. Let $C'_k= z_k+[\epsilon/8n, (1-\epsilon/8)/n]\times [\epsilon/n , (2-\epsilon)/n]$. Note that the sets $C'_k$ are disjoint, included in $\A_0 $, have the same $\Leb_{\A_0}$-volume equal to $(1-\epsilon)(1-\epsilon/4)/n^2$ and their diameters are $ < \sqrt5/n$. Likewise, for $v$ close enough to $(1-\epsilon)(1-\epsilon/4)/n^2$, there exists a disjoint family $(D'_k)_{1\le k\le N}$ such that each $D'_k\subset \A_0 $ is a neighborhood of $C'_k$ which is diffeomorphic to a closed disk and whose volume is $v$ and whose diameter is $<\sqrt 5/n$. 
 
By \cref{coro moser}, there exists $h \in \Sym^\infty(\A_0)$ which sends each $D_k$ to $D_k'$ and which is equal to the identity nearby $\partial \A_0$.  Thus $h\in \Ham^\infty(\A_0)$ and it can be extended to the identity to an element of $ \Ham^\infty_0(\A)$. 
 Then observe that the Kantorovitch-Wasserstein distance between 
$h_*\Leb_{\T\times \{y\}}$ and $N^{-1} \sum_{k=1}^N \delta_{z_k}$ is smaller than:
\[\frac1N \sum_{k=1}^N \sd(z_k, h(D_k) ) + \Leb_{\T\times \{y\}}(\T\times \{y\}\setminus \bigcup_k C_k) \le \max_k \diam D'_k + \epsilon/4 \le 
\sqrt 5 /n + \epsilon/4\]
which is smaller than $\epsilon/2$ for $n$ large enough. We recall that the distance between $N^{-1} \sum_{k=1}^N \delta_{z_k}$ and $\Leb_{\A_0}$ is $<\epsilon/2$ so the distance between $h_*\Leb_{\T\times \{y\}}$ and $\Leb_{\A_0}$ is $<\epsilon$. \end{proof}

\subsection{Lemma for emergence}\label{proof necklace}
\cref{neck lace} is:
 \begin{lemma} \label{neck lace2}
 For every $1>\epsilon>0$, there exist $\eta>0$ arbitrarily small, an integer $M\ge \exp( \eta^{-2+\epsilon})$, a map $h\in \Ham^\infty(\T\times \R/2\Z)$ and a family $( J_{\epsilon,i})_{1\le i\le M } $ of disjoint segments of $\R/2\Z$ such that:
\begin{enumerate}[(a)]
\item it holds $\leb(\I\setminus \bigcup_i J_{\epsilon,i} )= \frac{2}{M}\epsilon $ and $ \leb J_{\epsilon,i}= (1-\frac{\epsilon}{ M}) \frac{2 }{M} $ for every $i$,
 \item for any $j\neq i$, $y\in J_{\epsilon,i} $ and $y'\in J_{\epsilon,j}$, the distance between the measures $h_* \leb_{\T\times \{y\}} $ and $h_* \leb_{\T\times \{y'\}}$ is greater than $\eta$.
 \item $h$ coincides with the identity on neighborhoods of $\{0\}\times \R/2\Z$ and $\T\times \{0\}$. 
 \end{enumerate} 
 \end{lemma} 
\begin{proof}
The proof of this lemma combines \cref{coro moser} and a development of \cite[Prop 4.2]{BBochi21}. 
The present construction is depicted in \cref{figure_emergence}. We fix a positive $\epsilon<1/1000$. We will take a large integer $n\ge 3$ depending on $\epsilon$.
Let $N:= 2\cdot n^2$ and let $M$ be the integer part of $ (2N)^{-1/4} e^{ N/20}$: 
\[N:= 2\cdot n^2 \qand M= [ (2N)^{-1/4} e^{ N/20}]\; .\]
Now we consider the family of boxes $(B_{i,j})_{1\le i\le N/2, 1\le j\le M}$ defined by:
\[ B_{i,j} :=\left[ ({i+{\epsilon/2}})\frac2N, (i+1-{\epsilon/2})\frac2N \right]\times J_{\epsilon, j} \quad \text{with }
J_{\epsilon, j}:= \left[(j+\frac\epsilon{2M})\frac 2M, ({j+1-\frac{\epsilon}{2M}} ) \frac2M \right]\; .\]
The sets $( B_{i,j})_{i,j}$ are disjoint, included in $\T\times \R/2\Z \setminus (\{0\}\times \R/2\Z \cup \T\times \{0\})$ and have volume:
\[\leb B_\bullet :=\frac{2(1- {\epsilon} )}N \cdot \frac{2(1-\epsilon/{M})}M=\frac{2(1-{\epsilon} )}{n^2} \cdot \frac{(1-\epsilon/{M})}{M} .\] 
Note that $(J_{\epsilon,j})_j$ satisfies $(a)$. We are going to define a conservative map $h$ which sends each box $B_{i,j}$ into a union of $N$ much larger boxes $(\tilde B_k)_{1\le k\le N}$. 

Let us define $(\tilde B_k)_{1\le k\le N}$. As $N=2\cdot n^2$, there is a bijection:
\[k\in \{1,\dots, N\}\to (x_k,y_k)\in \{0,\dots, n-1\}\times \{0,\dots, 2n-1\}\; .\]
Let:
\[\tilde B_k= z_k+\left[\frac{\epsilon^2}2\frac1 n, \left(1-\frac{\epsilon^2 }2\right)\frac1 n \right]^2 \quad \text{with }z_k:= (\tfrac1n x_k\mod 1 , \tfrac1n y _k \mod 2)\in \T\times \R/2\Z\; .\]
The sets $(\tilde B_k)_k$ are $\epsilon^2/n$-distant, included in $\T\times \R/2\Z \setminus (\{0\}\times \R/2\Z \cup \T\times \{0\})$ and have volume:
\[\leb \tilde B_\bullet :=(1-\epsilon^2)^2/n^2.\] 
 
 Now we shall determine in which $\tilde B_k$ each $B_{i,j}$ is sent. To visualize this attribution, one can interpret each box $B_{i,j}$ as a colored pearl, each row $(B_{i,j})_{1\le i\le N/2}$ as a necklace and each $\tilde B_k$ as a package for pearls of the same color, see \cref{figure_emergence}. Hence a map $h$ which sends each $B_{i,j}$ into a certain $\tilde B_k$ define a coloring of each necklace. 
 We are going to fix a coloring given by the lemma below. In the language of colored necklaces, it states that $N$ colors 
 suffices to color $M$ necklaces of $N/2$ pearls each such that any pair of necklaces has at least a quarter of their pearls of different colors and such that each color is used approximately evenly among all the necklaces. 
\begin{lemma}\label{coloriage}For every $N=4\cdot n^2$ sufficiently large, there is a map:
\[C: \{1,\dots, N/2\}\times \{1,\dots, M\}\to \{1,\dots, N\}\]
 such that:
\begin{enumerate}[(1)]
\item For every $k\in \{1,\dots, N\}$, it holds $\# C^{-1}(\{k\}) \le M(1+{\epsilon/M}) /2 $. 
\item For every $j\neq j'$, it holds $\frac2N \#\{(i,i'): C(i,j)=C(i',j')\}<3/4$.\end{enumerate}
\end{lemma}
This lemma is a development of the second step of \cite[Prop 4.2]{BBochi21}. It is proved below. 

By (1), for every $1\le k\le N$, it holds:
\[\leb \bigcup_{(i,j)\in C^{-1}(\{k\}) } B_{i,j}= \# C^{-1}(\{k\}) \cdot \leb B_\bullet\le 
\frac M2 (1+{\epsilon/M}) \cdot \frac{2(1-{\epsilon} )}{n^2} \cdot \frac{(1-\epsilon/M)}M= \frac{ (1-\epsilon)(1-{\epsilon^2/M^2})}{n^2}\; .
\]
On the other hand, we have $\leb \tilde B_\bullet = (1-\epsilon^2)^2/n^2$.  
Hence, when $n$ is large ($\epsilon<0.001$ fixed), we have for every $k\in \{1,\dots, N\}$:
\[\leb \bigcup_{(i,j)\in C^{-1}(\{k\}) } B_{i,j}= \# C^{-1}(\{k\})\cdot \leb B_\bullet
<\leb \tilde B_k\; .\]
 So there exists $v$ greater and arbitrarily close to $\leb B_\bullet:=\frac{2(1-{\epsilon} )}{n^2} \cdot \frac{(1-\epsilon/{M})}{M}$ 
 such that for every $k$:
\[\# C^{-1}(\{k\}) \cdot v < \leb \tilde B_k = (1-\epsilon^2)^2/n^2\; .\]
For $v$ close enough to $\leb B_\bullet$, for every $(i,j)$, there is a neighborhood $D_{i,j}$ of $B_{i,j}$ of volume $v$ which is diffeomorphic to a closed disk, such that $(D_{i,j})_{i,j}$ are pairwise disjoint and disjoint from $\{0\}\times \R/2\Z \cup \T\times \{0\}$. 
Furthermore, for every $1\le k\le N$, there is a family $(\tilde D_{i,j})_{(i,j)\in C^{-1}(\{k\})}$ of disjoint disks which are included in $\tilde B_k$ and of volume $v$. 

Now we apply \cref{coro moser} which asserts the existence of $h\in \Ham^\infty(\T\times \R/2\Z)$ which sends each $D_{i,j}$ to $\tilde D_{i,j}$ and such that $h$ coincides with the identity on neighborhoods of $\{0\}\times \R/2\Z$ and $\T\times \{0\}$. Hence Property $(c)$ is satisfied. Let us verify Property $(b)$. 

For every $i$ and $y\in J_{\epsilon,j} $, the circle $\T\times \{y\}$ intersects each $D_{i,j}$ at a set of length $>2(1-{\epsilon})/N$ and no more than a subset of $\T\times \{y\}$ of measure $<\epsilon$ is not included in some $D_{i,j}$. 
Thus for any $j'\neq j$ and $y'\in J_{\epsilon,j'}$, 
the distance between the measures $h_* \leb_{\T\times \{y\}} $ and $h_* \leb_{\T\times \{y'\}}$ is at least:
\[ \min_{k\neq k'} d (\tilde B_k, \tilde B_{k'})\cdot \left(
\sum_{ i : C(i,j)\notin\{ C(i',j'): i'\} } \leb D_{i,j}\cap \T\times \{y\} - \epsilon\right)\]
\[\ge \frac{\epsilon^2}{n} \cdot \left( \frac 2N (1-\epsilon )\# \{ i : C(i,j)\neq C(i',j')\forall i'\}- 
\epsilon\right)
 \ge 
 \frac{\epsilon^2}{n} \cdot \left( \frac 2N (1-\epsilon )\frac N8- 
\epsilon\right)= \frac{\epsilon^2}{n} \cdot \left( \frac 14- \frac54 
\epsilon\right)\; .\]
As $\epsilon<1/1000$, the distance between the measures $h_* \leb_{\T\times \{y\}} $ and $h_* \leb_{\T\times \{y'\}}$ is at least:
 \[\eta:= \frac{\epsilon^2}{5\cdot n}, \]
as stated in $(b)$. It remains only to show the inequality $M\ge \exp( \eta^{-2+\epsilon})$. We recall that 
 \[M= [ (2N)^{-1/4} e^{ N/20}] = [ (2n)^{-1/2} e^{ n^2/10}] = \left[ \left( \frac{2\epsilon^2}{5\cdot \eta} \right)^{-1/2} e^{ \frac{\epsilon^4}{ 250\cdot \eta^2} }\right] \; .\]
 As $\epsilon$ is fixed, when $N$ is large, then $\eta$ is small and 
 so $ \frac{\epsilon^4}{ 250\cdot \eta^2}$ is large compared to $\tfrac12 \eta^{ -2+\epsilon}$ and $\sqrt{\frac{5\cdot \eta}{2\epsilon^2}}$ is large compared to $ \exp( \tfrac12 \eta^{ -2+\epsilon})$. 
 Thus $M$ is large compared to $ \exp( \eta^{-2+\epsilon})$.\end{proof}
\begin{proof}[Proof of \cref{coloriage}] We consider the space:
\[\cal F:= \{f: \{1,\dots, N\} \to \{0,1\} : \sum_{i=1} ^N f(i)=\tfrac N 2\}\; . \]
We can associate to each element $f\in \cal F$ a coloring of a necklace formed by $N/2$ pearls. Indeed we can set that the color $i$ is represented in the necklace iff $f(i)=1$. This defines an increasing map:
\[C_f: \{1,\dots, N/2\}\to \{1,\dots, N\} \; .\] 

We endow the space $\cal F$ with the Hamming distance:
\[\mathrm{Hamm} (f,g):= \#\{i \in \{1,\dots, N\}| f(i)\neq g(i)\}\; .\]
Now choose randomly $M$ elements $(f_j)_{1\le j\le M}$ of $\cal F$, for the equidistributed measure on the finite set $\cal F$ and we define:
\[C: (i,j)\in\{1,\dots, N/2\}\times \{1,\dots, M\}\mapsto C_{f_j} (i)\in \{1,\dots, N\}\; .\] 
Let us show that the map $C$ satisfies the statement of the lemma with positive probability.

 For every $k\in \{1,\dots, N\}$, the probability of the event $\{f_j(k)=1\}$ is $1/2$ for every $j$. Thus, with $\epsilon':= \epsilon/M$, by Bernstein's inequality the probability of the event $\{ \# C^{-1}(\{k\}) \le M(1+\epsilon') /2\} $ is at most: 
\[ p_0:= 2\cdot \exp \left(-{\frac {M\epsilon'^{2}}{2(1+{\frac {\epsilon'}{3}})}}\right).\]
Thus the probability that Property {\em (1)} fails is $\le N\cdot p_0\to 0$ when $N$ is large because $\log M\sim N/20$.

Given $j\neq j'$, the probability $p$ that $\#\{(i,i'): C(i,j)=C(i',j')\}\ge 3N/8$ is equal to $\# B/\#\cal F$ where $B$ is a ball of radius $N/4$ in $\cal F$. The following bound is obtained by first relaxing the condition $\sum f(i)=N/2$ and then applying Bernstein's inequality:
\begin{multline} 2^{-N}\cdot \# B \le \mathrm{Prob}\left\{ (X_n)_{1\le n\le N}\in \{0,1\}^{N}: \sum_n X_n \le \tfrac N4\right\}\\
= \mathrm{Prob}\left\{ (X_n)_{1\le n\le N}\in \{-1,1\}^{N}: \tfrac1N \sum_n X_n > \tfrac 12\right\}\le \exp\left(-\frac {N \cdot (1/2)^2}{2(1+\tfrac13\cdot \tfrac12 )}\right)\le \exp(-\tfrac 3{28N})\le \exp(- \tfrac N{10})\; . 
\end{multline}
On the other hand, the cardinality of $\cal F$ is $\binom{N}{N/2}\ge (2N)^{-1/2} 2^N$ (by Stirling's formula). Thus:
\[ p = \frac{\# B}{\#\cal F}\le \sqrt{2N}\cdot \exp(- \tfrac N{10})\; .\] 
The probability to having obtained a map $C$ which does not satisfy property {\em (2)} of the lemma is:
\[ \sum_{k=1}^{M} (k-1) \cdot p\le \tfrac12 M ^2 \cdot p\le \tfrac 12 \left( (2N)^{-1/4} e^{ N/20}\right) ^2\cdot 
\sqrt{2N}\cdot \exp(- \tfrac N{10}) =\tfrac12\; .\]
Hence the probability that both properties (1) and (2) fail is at most close to $1/2$. So there exists indeed $C$ satisfying both properties. 
\end{proof}
\bibliographystyle{alpha}
\bibliography{references.bib}

\def\cprime{$'$} \def\cprime{$'$} \def\cprime{$'$}
\begin{thebibliography}{BCLR06}

\bibitem[AK70]{AK70}
D.~Anosov and A.~Katok.
\newblock New examples in smooth ergodic theory, ergodic diffeomorphisms.
\newblock {\em Trans. Mosc. Math. Soc}, 23(1), 1970.

\bibitem[BB21]{BBochi21}
P.~Berger and J.~Bochi.
\newblock On emergence and complexity of ergodic decompositions.
\newblock {\em Advances in Mathematics}, 390:107904, 2021.

\bibitem[BB23]{BBiebler22}
P.~Berger and S.~Biebler.
\newblock Emergence of wandering stable components.
\newblock {\em J. Am. Math. Soc.}, 36(2):397--482, 2023.

\bibitem[BCLR06]{beguin2006pseudo}
F.~B{\'e}guin, S.~Crovisier, and F.~Le~Roux.
\newblock Pseudo-rotations of the open annulus.
\newblock {\em Bulletin of the Brazilian Mathematical Society}, 37(2):275--306,
  2006.

\bibitem[Bed18]{Be18}
E.~Bedford.
\newblock Fatou components for conservative holomorphic surface automorphisms.
\newblock In {\em Geometric Complex Analysis}, pages 33--54. Springer, 2018.

\bibitem[Ber17]{Be_steklov_17}
P.~Berger.
\newblock Emergence and non-typicality of the finiteness of the attractors in
  many topologies.
\newblock {\em Proceedings of the Steklov Institute of Mathematics},
  297(1):1--27, 2017.

\bibitem[Ber20]{BeJMP22}
P.~Berger.
\newblock Complexities of differentiable dynamical systems.
\newblock {\em Journal of Mathematical Physics}, 61(3):032702, 2020.

\bibitem[Ber23]{berger2021coexistence}
P.~Berger.
\newblock Coexistence of chaotic and elliptic behaviors among analytic,
  symplectic diffeomorphisms of any surface.
\newblock {\em J. {\'E}c. Polytech., Math.}, 10:525--547, 2023.

\bibitem[Ber24]{Be24}
P.~Berger.
\newblock Analytic pseudo-rotations {II}: a principle for spheres, disks and
  annuli.
\newblock {\em ArXiv}, 2024.

\bibitem[BGH25]{BGH22}
P.~Berger, N.~Gourmelon, and M.~Helfter.
\newblock Any diffeomorphism is a total renormalization of a close to identity
  map.
\newblock {\em Invent. Math.}, 239(2):431--468, 2025.

\bibitem[Bir41]{Bi41}
G.~Birkhoff.
\newblock Some unsolved problems of theoretical dynamics.
\newblock {\em Science}, 94(2452):598--600, 1941.

\bibitem[BK19]{BK19}
S.~Banerjee and P.~Kunde.
\newblock Real-analytic abc constructions on the torus.
\newblock {\em Ergodic Theory and Dynamical Systems}, 39(10):2643--2688, 2019.

\bibitem[BS06]{BSR2}
E.~Bedford and J.~Smillie.
\newblock Real polynomial diffeomorphisms with maximal entropy. {II}. {S}mall
  {J}acobian.
\newblock {\em Ergodic Theory Dynam. Systems}, 26(5):1259--1283, 2006.

\bibitem[BT24]{BT22}
P.~Berger and D.~Turaev.
\newblock {\em Israel Journal of Mathematics}, pages 1--16, 2024.

\bibitem[CP18]{CP18}
S.~Crovisier and E.~Pujals.
\newblock Strongly dissipative surface diffeomorphisms.
\newblock {\em Comment. Math. Helv.}, 93(2):377--400, 2018.

\bibitem[Cro06]{crovisier2006exotic}
S.~Crovisier.
\newblock Exotic rotations.
\newblock {\em Course Notes}, 2006.

\bibitem[Duj20]{Du20}
R.~Dujardin.
\newblock A closing lemma for polynomial automorphisms of {$\Bbb C^2$}.
\newblock Number 415, pages 35--43. 2020.
\newblock Some aspects of the theory of dynamical systems: a tribute to
  Jean-Christophe Yoccoz (volume I).

\bibitem[FH78]{FH78}
A.~Fathi and M.~R. Herman.
\newblock Existence de diff{\'e}omorphismes minimaux.
\newblock Ast{\'e}risque 49(1977), 37-59 (1978)., 1978.

\bibitem[FK04]{FK04}
B.~Fayad and A.~Katok.
\newblock Constructions in elliptic dynamics.
\newblock {\em Ergodic Theory and Dynamical Systems}, 24(5):1477--1520, 2004.

\bibitem[FK14]{FK14}
B.~Fayad and A.~Katok.
\newblock Analytic uniquely ergodic volume preserving maps on odd spheres.
\newblock {\em Commentarii Mathematici Helvetici}, 89(4):963--977, 2014.

\bibitem[Fur61]{Fu61}
H.~Furstenberg.
\newblock Strict ergodicity and transformation of the torus.
\newblock {\em American Journal of Mathematics}, 83(4):573--601, 1961.

\bibitem[Hel25]{He22}
M.~Helfter.
\newblock Scales.
\newblock {\em ArXiv}, 310(1):40, 2025.

\bibitem[Her98]{He98}
M.~Herman.
\newblock Some open problems in dynamical systems.
\newblock In {\em Proceedings of the International Congress of Mathematicians},
  volume~2, pages 797--808. Berlin, 1998.

\bibitem[Kat73]{Ka73}
A.~Katok.
\newblock Ergodic perturbations of degenerate integrable hamiltonian systems.
\newblock {\em Mathematics of the USSR-Izvestiya}, 7(3):535, 1973.

\bibitem[MS17]{mcduff2017introduction}
D.~McDuff and D.~Salamon.
\newblock {\em Introduction to symplectic topology}, volume~27.
\newblock Oxford University Press, 2017.

\bibitem[Nar71]{Na71}
R.~Narasimhan.
\newblock {\em Several complex variables}.
\newblock Chicago Lectures in Mathematics. University of Chicago Press,
  Chicago, Ill.-London, 1971.

\bibitem[Sul85]{sullivan1985quasiconformal}
D.~Sullivan.
\newblock Quasiconformal homeomorphisms and dynamics {I}. solution of the
  fatou-julia problem on wandering domains.
\newblock {\em Annals of mathematics}, 122(2):401--418, 1985.

\bibitem[Yoc84]{yoccoz1984n}
J.-C. Yoccoz.
\newblock Il n'y a pas de contre-exemple de {D}enjoy analytique.
\newblock {\em CR Acad. Sci. Paris S{\'e}r. I Math}, 298(7):141--144, 1984.

\end{thebibliography}

\end{document}